\documentclass[reqno,a4paper]{amsart}

\usepackage{enumerate}
\usepackage{amssymb}
\usepackage[all,ps,dvips]{xy} \SelectTips{cm}{}

\newtheorem{thm}{Theorem}[section]
\newtheorem{lem}[thm]{Lemma}

\newtheorem{prop}[thm]{Proposition}
\newtheorem{cor}[thm]{Corollary}

\newtheorem*{classthm}{Classification Theorem}

\theoremstyle{definition}
\newtheorem{rem}[thm]{Remark}

\newtheorem{example}[thm]{Example}

\newcommand{\Pic}{{\text{\rm Pic}}}
\newcommand{\Tors}{{\text{\rm Tors}}}

\newcommand{\V}{{\mathcal{V}}}
\newcommand{\OO}{{\mathcal {O}}}

\newcommand{\LL}{\mathcal{L}}

\newcommand{\epsi}{\epsilon}
\newcommand{\veps}{\varepsilon}
\newcommand{\fie}{\varphi}
\newcommand{\si}{\sigma}
\newcommand{\Si}{\Sigma}

\newcommand{\C}{{\mathbb C}}

\newcommand{\Q}{{\mathbb Q}}
\newcommand{\Z}{{\mathbb Z}}
\newcommand{\F}{{\mathbb F}}
\newcommand{\pp}{{\mathbb{P}}}

\newcommand{\DD}{{\mathcal{D}}}
\newcommand{\FF}{{\mathcal{F}}}

\newcommand{\T}{{\mathcal{T}}}

\newcommand{\tpi}{\tilde{\pi}}
\newcommand{\tsi}{\tilde{\si}}
\newcommand{\map}{\dasharrow}

\title{Numerical Godeaux surfaces with an involution}

\author[A.\ Calabri, C.\ Ciliberto and M.\ Mendes Lopes]{%
Alberto Calabri, Ciro Ciliberto and Margarida Mendes Lopes}

\thanks{
This research has been carried in the framework of the EAGER project
financed by the EC, project n.\ HPRN-CT-2000-00099.
The first two authors are members of G.N.S.A.G.A.-I.N.d.A.M.,
which generously supported this research.
The third author is a member of the Center for Mathematical Analysis,
Geometry and Dynamical Systems, IST  and was  partially  supported
by FCT (Portugal) through program  POCTI/FEDER and
Project POCTI/MAT/44068/2002.
}

\begin{document}

\begin{abstract}
Minimal algebraic surfaces of general type
with the smallest possible invariants
have geometric genus zero and $K^2=1$
and are usually called {\em numerical Godeaux surfaces}.
Although they have been studied by several authors,
their complete classification is not known.

In this paper we classify numerical Godeaux surfaces
with an involution, i.e.\ an automorphism of order 2.
We prove that they are birationally equivalent
either to double covers of Enriques surfaces,
or to double planes of two different types:
the branch curve either has degree 10
and suitable singularities, originally suggested by Campedelli,
or is the union of two lines and a curve of degree 12
with certain singularities.
The latter type of double planes are degenerations
of examples described by Du Val
and their existence was previously unknown;
we show some examples of this new type,
computing also their torsion group.
\end{abstract}

\maketitle

\vspace{-0.4cm}

\section{Introduction.}

In the one-century-and-a-half history of classification of algebraic varieties,
surfaces having geometric genus $p_g=0$
and irregularity $q=0$ have been studied from the very beginning.
They were conjectured to be rational by Max Noether (about 1870),
until Enriques, in 1894, suggested the existence of surfaces with $p_g=q=0$
and bi-genus $P_2=1$ which now bear his name.
After that, Castelnuovo, in 1896,
proved his celebrated \emph{rationality criterion}, which says
that a surface $X$ is rational if and only if $P_2(X)=q(X)=0$.
Since then, the classification of surfaces with $p_g=0$
has received particular attention by algebraic geometers, and not only.
Indeed, in the course of the years, it has been discovered
that these surfaces are interesting
not only for classification purposes, but also for their intriguing
relations with other fields of mathematics,
e.g.\ Bloch's conjecture, classification of four-folds,
etc.\ (see \cite{bloch, freedman}, for general information cf.\ \cite{survey}).

In 1931--32, Godeaux and Campedelli gave the first two examples
of minimal surfaces of general type with $p_g=0$.
Godeaux considered a quotient of a quintic surface in $\mathbb{P}^3$ by a
freely acting cyclic group of order 5 of projective transformations.
The smooth minimal model of this surface has $K^2=1$.
Campedelli constructed a double plane, i.e.\ a double cover of $\pp^2$,
branched along a degree 10 curve with six points, not lying on a conic,
all of type $[3,3]$,
that is a triple point with another infinitely near triple point.
Here the smooth minimal model has $K^2=2$.

Campedelli also proposed the construction of a minimal surface of
general type with $p_g=0$ and $K^2=1$ as the minimal model of a double plane branched
along a curve of degree 10 with a 4-tuple point and five
points of type $[3,3]$, not lying on a conic. The actual existence of
such a curve was proved only 50 years later by Kulikov, Oort and
Peters in \cite{OP}.
We will say that a double plane with branch curve having
the singularities suggested by Campedelli
is of \emph{Campedelli type}.

Minimal surfaces of general type with $p_g=0$ and $K^2=1$
are nowadays called \emph{numerical Godeaux surfaces}.
They have been studied classically, but also by several authors in the last 30 years:
it would be too long to recall here all the relevant contributions,
anyway most of them appear in our list of references.

As its construction shows, the original Godeaux's example
has non-trivial torsion, more precisely its torsion group
is cyclic of order 5.
It is actually a result of Miyaoka in \cite{miyaoka2}
that the torsion group of a numerical Godeaux surface is cyclic of order at most 5
and surfaces $S$ with $\Tors(S)=\Z/5\Z$ fill up an irreducible component
of the moduli space of the expected dimension $8=\chi(\T_S)$,
consisting of quotients of quintics in $\pp^3$ by a $\Z/5\Z$-action,
as in the original Godeaux's construction.

It has to be expected that the larger is the torsion,
the easier is the study and the classification of these surfaces.
This view-point has been pursued by Miles Reid in \cite{Re1},
who proved that also surfaces with $\Tors(S)=\Z/3\Z$ or $\Z/4\Z$
fill up an 8-dimensional irreducible component of the moduli space,
giving an explicit description of their canonical rings.

Although a few examples of numerical Godeaux surfaces $S$
with no torsion (cf.\ \cite{CG, DW}), or even simply connected
(see \cite{Ba2}), or with $\Tors(S)=\Z/2\Z$ (see \cite{Ba1, We1, We2}) are known,
neither a classification result,
nor a description of the moduli space are known in these cases.
Furthermore, all these examples turn out to possess an involution,
i.e.\ a birational automorphism of order 2.

This motivates the study of numerical Godeaux surfaces with an involution.
A first investigation of this subject has been done
by J.~Keum and Y.~Lee in \cite{keum},
under the assumption that the bicanonical system has no fixed components.
In this paper, we make no assumption of this sort and we prove the following:

\begin{classthm}
A numerical Godeaux surface $S$ with an involution $\sigma$
is birationally equivalent to one of the following:
\begin{enumerate}
\item a double plane of Campedelli type;
\item a double plane branched along
a reduced curve which is the union of two distinct lines $r_1,r_2$
and a curve of degree 12 with the following singularities:
\begin{itemize}
\item[$\bullet$] the point $q_0=r_1\cap r_2$ of multiplicity 4;
\item[$\bullet$] a point $q_i\in r_i$, $i=1,2$, of type $[4,4]$,
where the tangent line is $r_i$;
\item[$\bullet$] further three points $q_3,q_4,q_5$ of multiplicity 4
and a point $q_6$ of type $[3,3]$,
such that there is no conic through $q_1,\ldots,q_6$;
\end{itemize}
\item a double cover of an Enriques surface
branched along a curve of arithmetic genus 2.
\end{enumerate}
In case $(3)$, the torsion group of $S$ is $\Tors(S)=\Z/4\Z$,
whilst in case $(2)$ $\Tors(S)$ is either $\Z/2\Z$ or $\Z/4\Z$.
\end{classthm}

As we said, examples of surfaces of type (1) are known in the literature.
Surfaces of type (3) will be called of \emph{Enriques type}.
Examples of such surfaces have been produced by Keum and Naie
(cf.\ \cite{keum2, Naie}).
Double planes as in case (2) are, to the best of our knowledge, new in the literature.
They appear to be degenerations
of double planes with $p_g=4$ and $K^2=8$,
introduced by Du Val in 1952 (see \cite{DV}),
when classifying surfaces with non-birational bicanonical map
(cf.\ \cite{SantaCruz} for a modern reference).
For this reason, we call numerical Godeaux surfaces as in (2)
of {\em Du Val type}.
We give explicit examples of these surfaces, both with torsion group
$\Z/2\Z$ and $\Z/4\Z$, in section \ref{s:ex}.

In order to prove our classification theorem, we proceed as follows.
First of all, in section \ref{prima},
we prove some relevant properties of the fixed locus
of an involution acting
on a surface of general type with geometric genus zero.
In order to do so, we follow
ideas contained in joint work of the third author and Rita Pardini
(e.g.\ \cite{mp}),
namely we combine the topological and holomorphic fixed point formulas
with Kawamata-Viehweg's vanishing theorem.

Applying these results to numerical Godeaux surfaces $S$ with an involution,
in section \ref{sGodeaux} we give a rather precise
description of the fixed locus of the involution and
we prove that the bicanonical system is invariant under the involution.
This is actually the key ingredient for the proof of the classification
and explains why degenerations of Du Val double planes, which have non-birational
bicanonical map, come into play (cf.\ also recent results of Borrelli in \cite{borrelli}
on the classification of surfaces with non-birational bicanonical map and low invariants).

Since the pencil $|2K_S|$ is invariant, we can in fact consider its image
on the quotient surface under the involution,
which is a pencil $\DD$ of curves of arithmetic genus 2.
Using this, one sees that the quotient surface is either rational
or birational to an Enriques surface.
The latter case, i.e.\ the Enriques type,
is worked out in section \ref{enr},
where in particular we prove that numerical Godeaux surfaces of this type
have torsion group of order 4 and
are birational to the double cover of an Enriques surface branched
along a curve of arithmetic genus 2, with at most irrelevant singularities,
which moves in a linear system with no fixed component.

In case the quotient of $S$ by the involution is rational,
one studies the pencil $\DD$ using adjunction.
This leads to two different cases:
one in which the adjoint to this system is a base-point-free pencil of rational curves,
the other in which the adjoint is a pencil of curves of genus 1.
These two cases are analysed separately in sections \ref{s:camp} and \ref{duval},
respectively.
The former case leads to double planes of Campedelli type:
it suffices to suitably use the pencil of rational curves
to map the quotient surface to $\F_1$, and then to the plane.
The latter case leads to double planes of Du Val type.
Here the quotient surface is mapped to a weak Del Pezzo surface $X$,
i.e.\ $-K_X$ is big and nef,
having $K_X^2=1$ and four disjoint $(-2)$-curves, whose sum is an even divisor in $X$.
Rational surfaces with an even set of $(-2)$-curves have been studied
by Dolgachev, Mendes Lopes and Pardini in \cite{dmp} and, more recently,
by us in \cite{CCM}.
Indeed we apply the main result in \cite{CCM} to find a suitable birational morphism
of the weak Del Pezzo surface $X$ to $\pp^2$, which realizes
the original numerical Godeaux surface
as a Du Val double plane.

Then we
give more information about the previously unknown case
of numerical Godeaux surfaces of Du Val type.
In particular, in section \ref{s:torsion},
we examine the interplay between reducibility of the branch curve
and torsion.
Our main result in this direction is Theorem \ref{settica},
which, under some assumptions
on the branch curve, gives an useful criterion to decide whether
the torsion is $\Z/4\Z$ or $\Z/2\Z$, based on the existence or not
of plane curves of degree 8 with suitable singularities.

We remark that our results concern the birational classification
of \emph{pairs} $(S,\sigma)$, where $S$ is a numerical Godeaux surface
and $\sigma$ an involution of $S$.
We do not treat here, in general, the interesting problem of determining how many
involutions can occur on a given numerical Godeaux surface
and of which type according to our classification theorem.
However, our results do give some partial information.
For instance, if the ramification curve $R$ on $S$ has an irreducible
component of genus 2, then we are in case (3) (see Proposition \ref{p:exc}).
In any case, using the 2-torsion, we give
a criterion, i.e.\ Corollary \ref{c:notCamp},
based on the irreducibility of a certain plane cubic,
which allows us to distinguish between the Du Val and the Campedelli types.

As we said, in section \ref{s:ex}, we prove the existence of
numerical Godeaux double planes of Du Val type
and
we are able to compute the torsion of these examples.
We do this using Maple in two different ways.
One way is to find branch curves, with irreducible degree 12 component,
which are invariant under a projective
automorphism of order 2 of the plane,
which is an idea originally due to Stagnaro (cf.\ \cite{St}).
The resulting examples turn out to have torsion group $\Z/4\Z$.

Another way is to try and find the degree 12 component of the branch curve
suitably reducible in a line and an irreducible component of degree 11.
We find examples of this type with torsion $\Z/2\Z$ and with torsion $\Z/4\Z$
and
we prove that the former example does not have a different involution
which makes it a double plane of Campedelli type, therefore
it is certainly new in the literature.
Although we do not treat here moduli problems,
we prove that our examples both vary in families
whose images in the moduli space have dimension 5,
and the general member of each family is a Du Val double plane
with an irreducible degree 12 component of the branch curve
(cf.\ Corollaries \ref{l:5-dim} and \ref{l:irr}).

\section{Notation and conventions.} \label{notation}

In this section we fix the notation which will be used in this paper.

Let $S$ be a complex projective surface. We set:
\begin{tabbing}
\hspace{\parindent}\= $\Tors_n(S)$: \= explanation \kill
\> $\kappa(S)$: \> the Kodaira dimension of $S$; \\
\> $\chi(\FF)$: \> the Euler characteristic of a sheaf $\FF$ on $S$; \\
\> $\Pic(S)$: \> the Picard group of $S$; \\
\> $\Tors(S)$: \> the subgroup of $\Pic(S)$ composed of torsion elements; \\
\> $\Tors_n(S)$: \> the subgroup of $\Pic(S)$ composed of elements of torsion $n$; \\
\> $\rho(S)$: \> the rank of the N\a'{e}ron--Severi group of $S$; \\
\> $K_S$: \> a canonical divisor of $S$; \\
\> $p_g(S)$: \> the geometric genus of $S$, that is $h^0(S,\OO_S(K_S))$; \\
\> $q(S)$: \> the irregularity of $S$, that is $h^1(S,\OO_S)$; \\
\> $P_m(S)$: \> the $m$-th pluri-genus of $S$, that is $h^0(S,\OO_S(mK_S))$, $m\ge1$;
\end{tabbing}
If $S$ is clear from the context, sometimes we will write $p_g$, $q$, $K^2$, etc.,
instead of $p_g(S)$, $q(S)$, $K_S^2$, etc.

Let $X$ be a complex projective variety.
We denote by $e(X)$ the topological Euler characteristic of $X$
and by $p_a(X)$ the arithmetic genus of $X$.
Recall that if $D$ is a curve on a surface $S$,
then $p_a(D)=D(D+K_S)/2+1$.

We denote by $\equiv$ the linear equivalence of divisors on a surface
and by $\sim$ the numerical equivalence.
We usually omit the sign $\cdot$ of the intersection product of two divisors on a surface.

Recall that a \emph{$(-1)$-curve} is a smooth irreducible rational curve $C$
with $C^2=-1$.
More generally, ones says that a smooth irreducible rational curve $C$
with $C^2=-n<0$ is a \emph{$(-n)$-curve}.

We say that a divisor $D$ on a surface is \emph{nef and big} if $D^2>0$
and $D E\ge 0$ for every irreducible curve $E$.

If $x$ is a real number, we denote by $[x]$ its integer part,
i.e.\ the largest integer number less than or equal to $x$.

A singular point of type $[m,m]$ on a curve is a point of multiplicity $m$
with an infinitely near point again of multiplicity $m$.

\section{Involutions on surfaces.}\label{prima}

Let $S$ be a smooth, irreducible, projective surface over the
field $\C$ of complex numbers. An {\em involution} of $S$ is an
automorphism $\si$ of $S$ of order 2. Remark that if $S$ is a
minimal surface of general type, then any birational automorphism
is an isomorphism, therefore any birational automorphism of order
2 is an involution. If $X$ is any variety and $\psi:S\map X$ is a
rational map, one says that $\psi$ is \emph{composed} with the
involution $\si$ if $\psi\circ\si=\psi$.

Given an involution $\si$ on $S$, its fixed locus is the union of
a smooth, possibly reducible, curve $R$ and of $k$ isolated points
$p_1,\ldots,p_k$.

Let $\pi\colon S\to \Si:=S/\si$ be the quotient map. and set
$B:=\pi(R)$. The surface $\Si$ is normal and
$\pi(p_1),\ldots,\pi(p_k)$  are ordinary double points, which are
the only singularities of $\Si$. In particular, the singularities
of $\Si$ are canonical and the adjunction formula gives
$K_S\equiv\pi^*K_{\Si}+R$.

Let $\epsi\colon V\to S$ be the blowing-up of $S$ at
$p_1,\ldots,p_k$ and let $E_i$ be the exceptional curve over
$p_i$, $i=1,\ldots,k$. Then $\si$ induces an involution $\tsi$ of
$V$ whose fixed locus is the union of $R_0:=\epsi^*(R)$ and of
$E_1,\ldots,E_k$. Denote by $\tpi\colon V\to W:=V/\tsi$ the
projection onto the quotient and set $B_0:=\tpi(R_0)$,
$C_i:=\tpi(E_i)$, $i=1,\ldots,k$.
The surface $W$ is smooth and the $C_i$ are disjoint
$(-2)$-curves. Denote by $\eta\colon W\to \Si$ the map induced by
$\epsi$. The map $\eta$  is the minimal resolution of the
singularities of $\Si$ and there is a commutative diagram:
\begin{equation}\label{diagram}
\xymatrix{ V \ar[r]^{\epsi} \ar[d]_{\tpi} & S \ar[d]^{\pi} \\
W \ar[r]^{\eta} & \Si}
\end{equation}
The map $\tpi$ is a flat double cover branched on $\tilde
B=B_0+\sum_{i=1}^k C_i$, hence there exists a divisor $L$ on $W$
such that $2L\equiv \tilde B$ and
\begin{equation} \label{projf}
\tpi_*\OO_V=\OO_W\oplus \OO_W(-L).
\end{equation}

Also $K_V\equiv\tpi^*(K_W+L)$. With this notation:

\begin{prop}\label{kv}
Let $S$ be a minimal surface of general type and let $\si$ be an
involution of $S$. Then:
\begin{enumerate} [\rm(i)]
\item $2K_W+B_0$ is nef and big; \item $(2K_W+B_0)^2=2K_S^2$;
\item $H^i(W,\OO_W(2K_W+L))=0$, $i=1,2$.
\end{enumerate}
\end{prop}

\begin{proof}
By the adjunction formula and commutativity of diagram
(\ref{diagram}), we have
\begin{equation}\label{2K+B0}
\tpi^*(2K_W+B_0)\equiv 2K_V-2\sum_{i=1}^k E_i \equiv
\epsi^*(2K_S).
\end{equation}
Then $2K_W+B_0$ is nef and big because so is $2K_S$, proving part
(i). Statement (ii) also follows by formula \eqref{2K+B0}.

Finally for part (iii), we have the equivalence of $\Q-$divisors:
\[
K_W+L \equiv \frac{1}{2}(2K_W+B_0)+\frac{1}{2}\sum_{i=1}^k C_i .
\]
The divisor $(2K_W+B_0)/2=\eta^*(2K_{\Si}+B)/2$ is nef and big,
because so is $2K_W+B_0$, whereas $\frac{1}{2}\sum_{i=1}^k C_i  $
is effective, with zero integral part, and its support has normal
crossings. Thus $h^i(W,\OO_W(2K_W+L))=0$ for $i>0$ by
Kawamata--Viehweg vanishing theorem (see, e.g., Corollary 5.12,
c), of \cite{vieg}).
\end{proof}

For surfaces of general type with $p_g=0$ having an involution, one
can be more specific. As shown in \cite{dmp}, the holomorphic  and
topological  fixed point  formulas (see p.~566 in \cite{as}
and formula (30.9) in \cite{gre})
 yield:

\begin{lem}\label{numeri}
Let $S$ be a surface with $p_g=q=0$ and $\si$ an involution of
$S$. Then the number of isolated fixed points of $\si$ is
$k=K_SR+4$.
\qed
\end{lem}

Also the following properties, which will be very useful in the
sequel, hold.

\begin{prop}\label{formulas} Let $S$ be a minimal surface of general type
with $p_g=0$ and let $\si$ be an involution of $S$. Then:
\begin{enumerate} [\rm(i)]
\item $k\geq 4$; \item $K_W L+L^2=-2$; \item $h^0(W,
\OO_W(2K_W+L))=K_W^2+K_W L$; \item $K^2_W+K_W L\ge 0$;

\item $k=K_S^2+4-2h^0(W, \OO_W(2K_W+L))$.

\end{enumerate}
\end{prop}

\begin{proof}

Since $S$ is minimal of general type, $K_S$ is nef and so statement
(i) follows from Lemma \ref{numeri}.

 Since $S$ is of general type with $p_g(S)=0$, also
$q(S)=0$. Therefore:
\begin{equation}\label{pg-q}
p_g(\Si)=p_g(W)=0, \qquad q(\Si)=q(W)=0.
\end{equation}

By standard double cover formulas, we have
$\chi(\OO_V)=2\chi(\OO_W)+(L^2+K_W L)/2$, thus statement (ii)
follows from $p_g(W)=q(W)=0$.

By Proposition \ref{kv}, (iii), and the Riemann-Roch Theorem, one has:
\begin{align*}
h^0(W,\OO_W(2K_W+L))  &=  \chi(\OO_W(2K_W+L))
=(2K_W+L) (K_W+L)/2+1=  \\
  &=  K_W (K_W+L)+ L(K_W+L)/2+1.
\end{align*}
Thus (iii) follows by statement (ii), and (iv) is a trivial consequence of
(iii).

Finally, as  for statement (v),
it suffices to remember that $k=K_S^2-K_V^2$, $K_V^2=2(K_W+L)^2$
and use statements (ii), (iii).
\end{proof}

\begin{cor}\label{D}
Let $S$ be a minimal surface of general type with $p_g=0$ and let
$\si$ be an involution of $S$.  Then, with $W$ as above,
\begin{enumerate} [\rm(i)]
\item $|2K_W+B_0|\neq \emptyset$, $|2K_W+B_0+\sum_{i=1}^k
C_i|=|2K_W+B_0|+\sum_{i=1}^k C_i$ and
 $h^0(W,\OO_W(2K_W+B_0)=K_S^2+1-h^0(W,\OO_W(2K_W+L));$

\item if $D\in |2K_W+B_0|$, then $D$ is nef and big, 1-connected
and thus $h^0(D,\OO_D)=1$, $h^1(D,\OO_D)=p_a(D)$;
\item $p_a(D)>0$.
\end{enumerate}
\end{cor}

\begin{proof}
The Hurwitz formula gives
$ 2K_V \equiv \tilde{\pi}^*(2K_W+B_0+\sum_{i=1}^k C_i)$.
By the  projection formula \eqref{projf}, one has
\begin{equation}  \label{2K_V}
H^0(V,\OO_V(2K_V))= H^0(W,\OO_W(2K_W+L)) \oplus
H^0(W,\OO_W(2K_W+B_0+\sum_{i=1}^k C_i)).
\end{equation}
By Proposition
\ref{formulas}, (v) and  (i), $K_S^2-2h^0(W, \OO_W(2K_W+L))=k-4
\geq 0$.

Since $h^0(V,\OO_V(2K_V))= h^0(S,\OO_S( 2K_S))=K_S^2+1$, we
conclude that the second summand in \eqref{2K_V} is not 0.
For every $i=1, \ldots, k$,
one has that $C_i (2K_W+B_0+\sum_{i} C_i)=-2$,
and so
$
|2K_W+B_0+\sum_{i} C_i|=|2K_W+B_0|+\sum_{i} C_i.
$
This proves (i).

In Proposition \ref{kv}, (i), we already proved that $D$ is nef
and big. Therefore $D$ is 1-connected, see e.g.\ Lemma 2.6 in
\cite{Ma}, hence $h^0(D,\OO_D)=1$ and the final assertion of (ii)
follows by the Riemann--Roch Theorem.

Finally for assertion (iii), it suffices to note that,
by Propositions \ref{kv} and \ref{formulas}, $K_W D\geq 0$ and $D^2>0$.
\end{proof}

\begin{cor}\label{sic}
Let $S$ be a minimal surface of general type with $p_g=0$ and
$\si$ an involution of $S$. Then, with $W$ as above, $K_W^2\geq
K_V^2$.

\end{cor}

\begin{proof}
By Noether's formula,
because $\chi(\OO_W)=\chi(\OO_V)$
and $e(W)\le e(V)$, since
the map $\tilde \pi$ determines
an injection of $H^2(W,\C)$ into $H^2(V,\C)$.
\end{proof}

\begin{cor}\label{fixnumber}
Let $S$ be a minimal surface of general type with $p_g=0$, let
$\fie:S \to \pp^{K_S^2}$ be the \emph{bicanonical map} of $S$ and
let $\si$ be an involution of $S$. Then the following conditions
are equivalent:
\begin{enumerate} [\rm (i)]

\item  $\fie$ is composed with $\si$; \item
$h^0(W,\OO_W(2K_W+L))=0$; \item $K_W(K_W+L)=0$; \item
 the number of isolated fixed points of $\si$ is $k=K^2_S+4$.

\end{enumerate}
\end{cor}

\begin{proof}
The bicanonical $\fie$ is composed with $\si$ if and only if one of the summands in
formula \eqref{2K_V} above vanishes.
By Corollary \ref{D}, (i), the second summand is never 0 and so
$\fie$ is composed with $\si$ if  and only if
$h^0(W,\OO_W(2K_W+L))=0$, which in turn,  by Proposition
\ref{formulas}, (iii) and (v), is equivalent to $K_W(K_W+L)=0$ and
$k=K_S^2+4$.
\end{proof}

\begin{cor}\label{K2+4}
Let $S$ be a minimal surface of general type with $p_g=0$ and
$\si$ an involution of $S$. If the bicanonical map $\fie$ is
composed with $\si$, then:
\begin{enumerate} [\rm(i)]

\item $h^0(W,\OO_W(2K_W+B_0))=P_2(S)=1+K^2_S$;

\item for  $D\in |2K_W+B_0|$,  $h^0(D,\OO_D(D))=K^2_S$;

\item $K_W D=0$, $D^2=2K_S^2$ and $p_a(D)=K_S^2+1$;

\item $-4\le K_W^2\leq 0$, and $K_W^2=0$ if and only if $K_W\sim0$;

\item either $W$ is rational or the minimal model of $W$ is an
Enriques surface.
\end{enumerate}
\end{cor}

\begin{proof}

Part (i) is immediate from Corollaries \ref{D}, (i), and \ref{fixnumber}, (ii).
Part (ii) follows from part (i) by considering the long exact sequence
obtained from $$0\to \OO_W\to \OO_W(D)\to\OO_D(D)\to 0$$ and using
that $h^1(W,\OO_W)=0$. Now note that  $DK_W=2(K_W+L)K_W$,
so $K_WD=0$ follows from Corollary \ref{fixnumber}, (iii).
Note also that $D^2=2K_S^2$ by Proposition \ref{kv}, (ii),
therefore the adjunction formula implies that $p_a(D)=K_S^2+1$.
Since $k=K_S^2+4$, one has that $K_V^2=-4$
and so the first  inequality in (iv) follows from Corollary
\ref{sic}. The second inequality and the remainder of statement
(iv) follow from the Index Theorem because $D$ is nef and big by
Proposition \ref{kv}, (i).

Finally, again because $D$ is nef and big,  $K_W D=0$ implies that
$\kappa(W)\leq 0$. So statement (v) follows by the classification
of surfaces, since $p_g(W)=q(W)=0$.
\end{proof}

The following lemma will be used later.
We keep the notation introduced above.
Namely $D$ is a general member of the non-empty linear system $|2K_W+B_0|$.

\begin{lem}\label{E1}
Let $S$ be a minimal surface of general type with $p_g=0$ and
an involution $\si$ and let $W$ be as above. If $E\subset W$ is a  curve such that $E
D=0$, then $E^2<0$ and the intersection form on the
components of $E$ is negative definite. In particular if $E$ is a
curve such that $E^2=-1$ and $E D=0$, then  $E
(\sum_{i=1}^k C_i)\leq  0$.
\end{lem}

\begin{proof}
The first part of the lemma is obvious by the Index Theorem,
because $D$ is nef and big by Proposition \ref{kv}, (i).
For the second part, note that, since $E(B_0+\sum_{i=1}^k C_i)$ is even and
$ED=0$, also $E (\sum_{i=1}^k C_i)  $ is even. By the
first part applied to the curve $E+C_i$,
$i=1,\ldots,k$, one has $E C_i\leq 1$. Suppose that $E
(\sum_{i=1}^k C_i)>0$. Then there would be at least two curves, say
$C_1$ and $C_2$, such that $EC_1=EC_2=1$, therefore $A=2E+C_1+C_2$
would satisfy $A^2=0$ and $AD=0$, a contradiction.
\end{proof}

\begin{prop}\label{W'}
Let $S$ be a minimal surface of general type with $p_g=0$ and
an involution $\si$  and let $W$ be as above.
Then there exists a birational morphism $f:W\to
W'$ and an effective divisor $D'$ on $W'$ with the following properties:
\begin{enumerate} [\rm (i)]
\item there are $k$ $(-2)$-curves $C'_i$ on $W'$ such that
$f^*(C'_i)=C_i$, $i=1, \ldots, k$;
\item $D'$ is nef
and such that $f^*(D')=D$, $D'^2=D^2$ and $p_a(D')=p_a(D)$;
\item
$K_{W'}+D'$ is nef.
\end{enumerate}
\end{prop}

\begin{proof}
If $K_W+D$ is nef, there is nothing to prove. Otherwise, since
$K_W+D$ is effective (see Corollary \ref{D}), there is an irreducible curve $E_1$ such
that $E_1^2<0$ and $E_1(K_W+D)<0$. Then $K_W E_1<0$,
thus $E_1$ is a $(-1)$-curve and $E_1 D=0$. Hence
$B_0 E_1=2$, and furthermore $E_1\cap \bigcup_{i=1}^k C_i =\emptyset$
by Lemma \ref{E1}.

Let $f_1:W\to W_1$ be the contraction of $E_1$ to a point $p_1$.
On $W_1$ we have $k$ $(-2)$-curves $C^{(1)}_i$ such that
$f_1^*(C^{(1)}_i)=C_i$, $i=1 \ldots k$, and a nef divisor $D_1$ on
it such that $D_1^2=D^2$ and $f_1^*(D_1)=D$. We set
$B^{(1)}_0=f_1(B_0)$.
Note that $B^{(1)}_0$ has a double point at
$p_1$ and it does not meet $\cup_{i=1}^k C^{(1)}_i$. Note also
that $D_1\equiv 2K_{W_1}+B^{(1)}_0$.

If $K_{W_1}+D_1$ is nef, we have finished. Otherwise, reasoning as
above, there is a $(-1)$-curve $E_2$ on $W_1$ such that
$E_2(K_{W_1}+D_1)<0$, $D_1 E_2=0$ and $E_2B_0^{(1)}=2$.
The curve $E'= f_1^*(E_2)$ satisfies
 $E' D=0$, $E'^2=-1$
and $K_W E'= -1$ and, since $E_1$ is disjoint from the curve $\sum_{i=1}^k C^{(1)}_i$,
none of the curves $C_i$ is a component of $E'$.
Hence, by Lemma \ref{E1}, $E' \cap \sum_{i=1}^k C^{(1)}_i=\emptyset$, implying that also $E_2$ is
disjoint from the curves $C^{(1)}_i$.
So we can contract $E_2$ and proceed.

Finally, the existence of the morphism $f:W\to W'$ is shown by
iterating the above procedure.
\end{proof}

\begin{rem} \label{divisibilita} We notice that, by  the proof of Proposition \ref{W'}, the
 divisor $D'+ C'_1+\cdots+C'_k$
  is divisible by $2$ in $\Pic(W')$.
\end{rem}

Later, we will need the following lemma, due to
Beauville (Lemme 2 in \cite{beauville}):

\begin{lem} \label{l:2-tors}
Let $\tilde\pi:V\to W$ be a flat double cover between two smooth
surfaces $V$ and $W$, branched over the smooth curve $\tilde
B=B_1+\cdots+B_n$, where $B_1,\ldots,B_n$ are the irreducible
components of $\tilde B$. Suppose that $\Pic(W)$ has no 2-torsion
element. Define a group homomorphism
$\psi: \Z_2^{\oplus n} \to \Pic(W)\otimes\Z_2$
by $\psi(\veps_1,\ldots,\veps_n) = \sum_{i=1}^n \veps_i B_i$.
Then $\Tors_2(V) \cong \ker(\psi)/\langle (1,\ldots,1) \rangle$.
\hfill$\Box$
\end{lem}

\section{Numerical Godeaux surfaces with an involution.}\label{sGodeaux}

In the remainder of this paper, we focus on the study of a \emph{numerical
Godeaux} surface $S$, i.e.\ a minimal surface of general type
with $p_g=0$ and $K_S^2=1$, having an involution $\si$.
We will freely use the notation introduced so far.

We start by recalling the following:

\begin{thm}[{Miyaoka, \cite{miyaoka2}}] \label{t:tors}
Let $S$ be a numerical Godeaux surface.
Then the torsion group $\Tors(S)$ of $S$ is cyclic of order $n\le 5$.
Moreover the linear system $|3K_S|$ has no fixed part and has at most 2 base points.
If there is no base point, then $|\Tors(S)|\le 2$.
If there is one base point, then $3\le |\Tors(S)|\le 4$.
If there are two base points, then $|\Tors(S)|=5$.
\qed
\end{thm}

\begin{lem}[{Lemma 5 in \cite{miyaoka2}}]\label{lemma5}
Let $S$ be a numerical Godeaux surface.
If $D$ is an effective divisor with $h^0(S,D)\geq2$, then $DK_S\ge2$.
\qed
\end{lem}

\begin{lem}[{Reid, see p.\ 158 in \cite{D}}] \label{d:tors}
Let $S$ be a numerical Godeaux surface.
If $\eta$ is a non-trivial element of $\Tors(S)$,
then there is a unique element $P_\eta$ in $|K_S+\eta|$.
Furthermore, if $\eta$ and $\eta'$, $\eta\neq \eta'$, are
non-trivial  elements of  $\Tors(S)$,
then $P_{\eta}$ and $P_{\eta'}$ have no common components.
\qed
\end{lem}

As an immediate consequence of Lemma \ref{d:tors}, we have the following:

\begin{cor}\label{components}
Let $S$ be a numerical Godeaux surface.
If $|\Tors(S)|= 4$ or $5$, then $|2K_S|$ has no fixed components.
\qed
\end{cor}

Our results in \S\ref{prima} imply that $S$ enjoys the
following properties (cf.\ also \cite{keum}):

\begin{prop}\label{godeaux}
If $S$ is a numerical Godeaux surface with an involution $\si$,
then:
\begin{enumerate} [\rm(i)]
\item the number of isolated fixed points of $\si$ is $k=5$;
\item the bicanonical map $\fie$ is composed with $\si$;
\item $K_S R=1$;
\item $R^2$ is odd and $-7 \leq R^2 \leq 1$.
\end{enumerate}
Furthermore $R=\Gamma+Z_1+\cdots+Z_h$ where:
\begin{enumerate} [\rm(i)] \addtocounter{enumi}{4}
\item $\Gamma$ is a smooth curve with $K_S \Gamma=1$, $0\le p_a(\Gamma)\le 2$
and $\Gamma^2=2p_a(\Gamma)-3$;
\item if $p_a(\Gamma)=2$, then $\Gamma\sim K_S$ and $S$ has non-trivial torsion;
\item $Z_1,\ldots,Z_h$ are disjoint $(-2)$-curves,
which are disjoint also from $\Gamma$, and
\begin{equation} \label{e:h}
h=p_a(\Gamma)-K_W^2-2 \ge 0.
\end{equation}

\end{enumerate}
\end{prop}

\begin{proof}
By Proposition \ref{formulas}, (i) and (v), the number $k$ is odd
and $4\leq k\leq 5$.
Hence $k=5$, that is part (i). Then, part  (ii) follows
by Corollary \ref{fixnumber},
whilst part (iii) follows from (i) and Lemma \ref{numeri}.

Let us prove part (iv). Since $2R^2=B^2$, one has $4L^2=2R^2-10$
and thus $L^2=(R^2-5)/2$.
By Proposition \ref{formulas}, (ii) and   Corollary \ref{fixnumber}, (iii),
$L^2=-2-K_W L=-2+K_W^2$, hence $R^2=1+2K_W^2$.  Then  by Corollary \ref{K2+4}, (iv), $-7\leq R^2\leq 1$.

Part (iii) implies that $R$ has a unique irreducible component
$\Gamma$ such that $K_S\Gamma=1$.
Since $R$ is smooth, so is $\Gamma$.
We can write $R=\Gamma+Z$, where $Z$ is effective.
Then $K_S Z=0$, thus the irreducible components of $Z$ are $(-2)$-curves,
which are pairwise disjoint and disjoint from $\Gamma$ because $R$ is smooth.

By the Index Theorem, $\Gamma^2\leq 1$.
Since $K_S\Gamma=1$ by adjunction, one has $\Gamma^2=2p_a(\Gamma)-3$,
which ends the proof of statement (v).

If $\Gamma^2=1$, then $\Gamma$ is homologous to $K_S$, but not
linearly equivalent to $K_S$, because $p_g=0$. Therefore
$\Gamma-K_S$ is a non-trivial torsion element of the
N{\'e}ron-Severi group of $S$, which shows (vi).

Finally (vii) follows from $R^2=1+2K_W^2=\Gamma^2-2h$
and $\Gamma^2=2p_a(\Gamma)-3$.
\end{proof}

\begin{rem}\label{r:godeaux}
The statements about $R$ on $S$ in Proposition \ref{godeaux}
can be read as well as about the curve $B_0$ on
the surface $W$,
cf.\ diagram \eqref{diagram} and Corollary \ref{K2+4}.
For instance, $\tilde\pi(\epsi^{-1}(Z_i))$ is
a smooth rational curve on $W$ with self-intersection $-4$
and is an irreducible component of $B_0$.
\end{rem}

Now Beauville's Lemma \ref{l:2-tors}
and Theorem \ref{t:tors} imply the following:

\begin{cor}\label{c:hgodeaux}
Let $h$ be the number of $(-2)$-curves of $R$ as in Proposition \ref{godeaux}.
If  $W$ is a rational surface, then
\begin{equation} \label{e:h2}
h \le 1 +\left[ \frac{-K_W^2}{2} \right]
\end{equation}
and if equality holds, then $S$ has non-trivial 2-torsion.
\end{cor}

\begin{proof}
Consider the map $\psi$ of Lemma \ref{l:2-tors}.
The domain of $\psi$ is $\Z_2^{\oplus(h+6)}$, because
the branch locus of $\tpi:V\to W$ has $h+6$ irreducible components.
Moreover the image of $\psi$ is a totally isotropic subspace
of $\Pic(W)\otimes\Z_2=H^2(W,\Z_2)$.

Since $W$ is simply connected, $h^2(W,\Z_2)=h^2(W,\Z)$.
Noether's formula implies that $b_2(W)=10-K_W^2$, hence
a totally isotropic subspace of $H^2(W,\Z_2)$
has dimension at most $5+[-K_W^2/2]$.
By Theorem \ref{t:tors}, $\dim\ker(\psi)\le 2$,
therefore $h+6\le 2+5+[-K_W^2/2]$, which proves \eqref{e:h2}.
If equality holds in \eqref{e:h2}, then
$\dim\ker(\psi)=2$, so $\Tors_2(V)$ and $\Tors_2(S)$ are not trivial
by  Lemma \ref{l:2-tors}.
\end{proof}

In order to prove the Classification Theorem,
stated in the introduction, we need to understand the surface $W$ and the divisor
$D\equiv 2K_W+B_0$ on $W$. By Proposition \ref{godeaux}
and Corollaries \ref{D} and \ref{K2+4},
one has that $D$ is nef, $D^2=2$, $K_W D=0$, $p_a(D)=2$, and
$h^0(W,\OO_W(D))=2$.

If $R$ has a component $\Gamma$ with $p_a(\Gamma)=2$,
we have further information on $W, D$:

\begin{cor}\label{c:godeaux}
Let $S$ be a numerical Godeaux surface with an involution $\si$.
Suppose that $R$ has an irreducible component $\Gamma$ of genus 2.
Then either one of the following two cases occurs:
\begin{enumerate}[\rm (i)]
\item $W$ is a minimal Enriques surface, $R=\Gamma$ and $|D|=|B_0|$;
\item $W$ is a rational surface, $-2\le K_W^2\le -1$ and
$\Gamma-K_S\in\Tors_2(S)$.
\end{enumerate}

\end{cor}

\begin{proof}
On $W$, let $\Gamma_0=\eta^*(\pi(\Gamma))\leq B_0$,
thus $\Gamma_0^2=2$ and $K_W\Gamma_0=0$.
Write:
\begin{equation}\label{Gamma_0}
D\equiv 2K_W+B_0 \equiv 2K_W+\Gamma_0+(B_0-\Gamma_0).
\end{equation}
Now $\Gamma_0 D = \Gamma_0^2=D^2$,
so, by the Index Theorem, $\Gamma_0 \sim D$.
Hence \eqref{Gamma_0} implies that $2K_W+B_0-\Gamma_0 \sim 0$.
By Corollary \ref{K2+4}, (v), either the minimal model of $W$ is an Enriques surface
or $W$ is rational.
In the former case $2K_W$ is an effective divisor.
Since $B_0-\Gamma_0$ is also effective,
$2K_W+B_0-\Gamma_0 \sim 0$ implies that  $0 \equiv 2K_W \equiv B_0-\Gamma_0$,
that is case (i) of the statement.

If $W$ is rational, numerical equivalence is the same
as linear equivalence and so $\Gamma_0 \equiv D$ and $B_0-\Gamma_0 \equiv -2K_W$,
therefore $2\Gamma \equiv 2K_S$.
By formula \eqref{e:h}, $B_0$ has $1-K_W^2$ irreducible components,
i.e.\ $h=-K_W^2$,
thus formula \eqref{e:h2} implies $K_W^2\ge -2$.
\end{proof}

\begin{rem}
We will see later, in Proposition \ref{p:exc}, that case (ii) in
Corollary \ref{c:godeaux} does not actually occur.
\end{rem}

We now study $W$ and
the pencil $|D|$ on $W$ by using adjunction:

\begin{lem}\label{emme}
With the above notation,
write $|K_W+D|=F+|M|$, where $F$ denotes the
fixed part and $|M|$  the movable part.
Then:
\begin{enumerate}[\rm (i)]
\item  $h^0(W,\OO_W(M))=2$;
\item  the general curve of the pencil $|M|$ is irreducible;
\item  $M D=2$;
\item  if $F\neq 0$, then every component $E$ of $F$ is
such that $D E=0$ and $E^2<0$.
\end{enumerate}
\end{lem}

\begin{proof}
Assertion (i) follows since $p_a(D)=2$ and $W$ is regular.
Assertion (ii) follows by Bertini's Theorem, since every
pencil on the regular surface $W$ is rational.

Let us prove part (iii).
Since $D$ is nef and $D (M+F)=2$, one has $M D\leq 2$.
Suppose by contradiction that $MD\leq 1$. Consider the pull--back
$|\tilde M|$ of $|M|$ to $V$, which is also the pull--back of a
pencil $|N|$ on $S$. Since $\tilde\pi^*(D)=\epsi^*(2K_S)$, one
would have that $N K_S \leq 1$, which is impossible by Lemma \ref{lemma5}.
This proves (iii).

Assertion (iv) follows now by the Index Theorem.
\end{proof}

We have already seen that $-4\leq K_W^2\leq 0$.
Now we want to consider the surface $W'$ as in Proposition \ref{W'}.

\begin{lem}\label{K2W'}
One has $-2\le K_{W'}^2 \le 0$. Furthermore
$K_{W'}^2=0$ if and only if $W'$ is a minimal Enriques surface.
\end{lem}

\begin{proof}
Since $D' (K_{W'}+D')=2$, the Index Theorem implies that
$(K_{W'}+D')^2\leq 2$, or equivalently $K_{W'}^2\le 0$, because
$K_{W'} D'=0$ and $D'^2=2$. On the other hand,
$(K_{W'}+D')^2\ge 0$ and therefore $K_{W'}^2\ge -2$.

If $K_{W'}^2=0$, then $K_{W'}\sim 0$, and therefore $W'$ is an
Enriques surface.
The converse is trivial.
\end{proof}

\begin{lem} \label{partfix}
If $K_{W'}^2<0$, then $|K_{W'}+D'|=|M'|$ has no fixed part.
\end{lem}

\begin{proof}
Write, as usual, $|K_{W'}+D'|=F'+|M'|$, where $F'$ is the fixed
part and $|M'|$ is the movable part. By Lemma \ref{emme}, (iii), and the
construction of the morphism $f:W\to W'$, which contracts only
curves in $F$, we see that $M' D'=2$.

Notice that $F' K_{W'}=F' M'+F'^2$, so $M'
F'$ is even. Since $M' F'\leq M' F'+M'^2\leq
(K_{W'}+D')^2=K_{W'}^2+2<2$, it follows that $M' F'=0$. This
forces $F'$ to be $0$. Otherwise, since $D' F'=0$, then
$F'^2<0$, implying that $F' (K_{W'}+D')= F'^2+ M'
F'<0$, a contradiction because $(K_{W'}+D')$ is nef.
\end{proof}

\begin{rem}\label{KW>-4}
Since $K_W^2\ge-4$ by Corollary \ref{K2+4}, (iv), the birational map $f:W\to W'$ of Proposition \ref{W'}
is the contraction of at most 3 (resp., at most 4) exceptional curves if
$W'$ is rational (resp., an Enriques surface).

Note that an unessential singularity of $B'_0$ corresponding
to a triple point requires at least 4 blowing-ups to be resolved,
and exactly 4 unless there is a double point infinitely near
to the triple point.
Thus $B'_0$ can have a triple point only if $W'$ is a minimal
Enriques surface, $K_W^2=-4$ and $B'_0$ has only one triple point
with no infinitely near double point.

If $K_W^2=-4$, then $\tilde\pi$ induces an isomorphism between
$H^2(W,\Z)$ and $H^2(V,\Z)$.
This implies that there is no rational curve on $W$ that does not meet
the branch locus $\tilde B$.
As a consequence, an unessential singularity of $B'_0$
can only be an ordinary triple point, a node or a cusp.
\end{rem}

\section{On numerical Godeaux surface of Enriques type.}\label{enr}

In this section we keep the notation we
introduced in the above sections.

Let $S$ be a numerical Godeaux surface with an involution.
We analyse here the case in which $W$ is birational to an
Enriques surface (see Corollary \ref {K2+4}, (v)).
This situation corresponds to the case in which $W'$ is a minimal
Enriques surface and $B'_0=D'$ (see Lemma \ref{K2W'}). The main
information is given by the following result:

\begin{prop}\label{enr1}
In the above setting, the general curve in the linear system $|D'|$ is
smooth and irreducible.
\end{prop}

\begin{proof} Since $D'^ 2=2$ and $h^ 0(W',\mathcal O_{W'}(D'))=2$
by Lemma \ref{emme} and Proposition \ref{W'},
it suffices to prove that the linear system $|D'|$ has no fixed component.
Assume the contrary and write $|D'|=\Phi+|\Psi|$, where $\Phi$ is the fixed part.
Since $D'$ is nef, then $D'\Psi\leq 2$.
On the other hand, by Lemma \ref{lemma5}
we have $D'\Psi=2$. This
implies $\Psi^2=0$, $\Psi \Phi=2$, $\Phi^2=-2$. Moreover by Proposition \ref{W'}, one
has $B'_0\cap (C'_1+\dots+C'_5)=\emptyset$, which implies
$\Psi(C'_1+\dots+C'_5)=0$. Now we remark that $|\Psi|$ is a pencil of elliptic
curves on the Enriques surface $W'$,
thus $\Psi$ is divisible by $2$ in $\Pic(W')$.
Since $\Psi+\Phi+C'_1+\dots+C'_5$ is
also divisible by $2$, we have that $\Phi+C'_1+\dots+C'_5$ is divisible by $2$.
Since $\Psi(\Phi+C'_1+\dots+C'_5)=\Psi\Phi=2$, we find a contradiction.
\end{proof}

\begin{rem}\label{enr+}
The above Proposition tells us also that every numerical Godeaux surface with
an involution of Enriques type is deformation equivalent to the double
cover of a minimal Enriques surface $W$ with five nodes $C_1,\dots,C_5$,
branched along a  smooth curve $B_0+C_1+\dots+C_5$, where $B_0$ is irreducible
of genus $2$. There are examples of surfaces of this type
(see Example 4.3 of \cite{keum}, cf.\ Corollary \ref{c:godeaux}
and Proposition \ref{p:exc} below).
\end{rem}

Next we describe the torsion group of these surfaces.

\begin{prop}\label{enrtors}
The torsion group of a numerical Godeaux surface $S$
of Enriques type is $\Z_4$.
\end{prop}

\begin{proof}
By Remark \ref{enr+}, we may assume that $W$
is a minimal Enriques surface and that $B_0$ is a smooth, irreducible curve
of genus $2$.
One has a cartesian diagram:
\begin{equation}\label {enrdiagram}
\xymatrix{%
Y \ar[r] \ar[d]_p & V \ar[d]^{\tilde\pi} \\
T \ar[r]^g & W
}
\end{equation}
where $g: T\to W$ is the K3--double cover of $W$, so $Y\to V$
is an \'etale double cover, where the minimal model of $Y$
has $K^ 2=2$ and $p_g=1$.
Thus $V$, and hence $S$, has $2$--torsion.
Moreover $\Tors(V)=\Z_4$ if and only if
$\Tors(Y)=\Z_2$.

Let us now look at the double cover $p: Y\to T$
which is branched along  the curve $H:=g^ *(B_0)$ and the ten $(-2)$--curves
which are pull-back via $g$ of the five $(-2)$--curves $C_1,\dots,C_5$.
Standard double cover considerations show that the bicanonical map of $Y$
factors through $p$, in particular it is not birational.
By Theorem 6.1 of \cite{CaDe}, this
implies that $Y$ has $2$--torsion, which implies the assertion.
\end{proof}

\begin{rem}
By Proposition \ref{enrtors} and Corollary \ref{components},
it follows that, if $S$ is a numerical Godeaux surface
of Enriques type, then $|2K_S|$ has no fixed component.
\end{rem}

\section{On numerical Godeaux surfaces of Campedelli type.}\label{s:camp}

In this section we follow the notation introduced above and we
study the case $K_{W'}^2=-2$ (cf.\ Lemma \ref{K2W'}).
As we will see, this case
corresponds to numerical Godeaux surfaces which are birationally
equivalent to double planes branched along a curve of degree 10
with a point of multiplicity 4 and five points of type $[3,3]$. As
we said in the introduction, we call these surfaces
of \emph{Campedelli type}.

\begin{lem} \label{l:M'rational}
Let $S$ be a numerical Godeaux surface with an involution
and let $W'$ as in Sections \ref{prima}, \ref{sGodeaux}.
Suppose that $K_{W'}^2=-2$.
Then $|M'|$ is a base
point free pencil of rational curves. Furthermore, if $A$ is an
irreducible component of a reducible fibre of $|M'|$, then one of
the following occurs:
\begin{enumerate} [\rm(i)]
\item either $A$ is a $(-1)$-curve such that $D' A=1$ and $B'_0
A=3$; \item or $A$ is a $(-2)$-curve such that $D' A=B'_0 A=0$.
\end{enumerate}
\end{lem}

\begin{proof}
The hypothesis $K_{W'}^2=-2$ implies that $M'^2=0$ and $K_{W'}
M'=-2$, so that $|M'|$ is a base point free pencil of rational curves.

An irreducible component of a reducible fibre of $|M'|$ is a
rational curve $A$ with $A^2 < 0$. Statements (i) and (ii) follow
from $0=M' A=(K_{W'}+D') A \ge K_{W'} A=-2-A^2$, where the
inequality holds because $D'$ is nef.
\end{proof}

Lemma \ref{l:M'rational} immediately implies the following corollary
(cf.\ also \cite{xiao}, \cite{horikawa}).

\begin{figure}[ht]
\vspace{-5mm}
\begin{align*}
& & (i) \qquad &
\xymatrix@!0@C+10pt{*=0{\circ}
\ar@{-}[]+<1.6pt,0pt>;[rr]-<1.6pt,0pt> & *=0{\bullet} & *=0{\circ}
\ar@{-}[]+<1.6pt,0pt>;[rr]-<1.6pt,0pt> & *=0{\bullet} & *=0{\circ}
\ar@{-}[]+<1.6pt,0pt>;[rr]-<1.6pt,0pt> & *=0{\bullet} & *=0{\circ}
\save"1,1"+<0ex,2ex>*{-1} \restore
\save"1,2"+<0ex,2ex>*{-2} \restore
\save"1,3"+<0ex,2ex>*{-2} \restore
\save"1,4"+<0ex,2ex>*{\cdots} \restore
\save"1,5"+<0ex,2ex>*{-2} \restore
\save"1,6"+<0ex,2ex>*{-2} \restore
\save"1,7"+<0ex,2ex>*{-1} \restore
\save"1,1"-<0ex,2ex>*{1} \restore
\save"1,2"-<0ex,2ex>*{1} \restore
\save"1,3"-<0ex,2ex>*{1} \restore
\save"1,4"-<0ex,2ex>*{\cdots} \restore
\save"1,5"-<0ex,2ex>*{1} \restore
\save"1,6"-<0ex,2ex>*{1} \restore
\save"1,7"-<0ex,2ex>*{1} \restore }
\\[-3mm]
& & \raisebox{-16pt}{$(ii)$} \qquad &
\xymatrix@!0@C+10pt@R-8pt{ & & & & & & & *=0{\circ} \\
*=0{\circ} \ar@{-}[]+<1.6pt,0pt>;[rr]-<1.6pt,0pt>
 & *=0{\bullet} & *=0{\circ} \ar@{-}[]+<1.6pt,0pt>;[rr]-<1.6pt,0pt>
 & *=0{\bullet} & *=0{\circ} \ar@{-}[]+<1.6pt,0pt>;[rr]-<1.6pt,0pt>
 & *=0{\bullet} & *=0{\circ}
  \ar@{-}[]+<1.5pt,0.9pt>;[ur]-<1.5pt,0.9pt>
  \ar@{-}[]+<1.5pt,-0.9pt>;[dr] & \\
 & & & & & & & *=0{\bullet}
\save"2,1"+<0ex,2ex>*{-1} \restore
\save"2,2"+<0ex,2ex>*{-2} \restore
\save"2,3"+<0ex,2ex>*{-2} \restore
\save"2,4"+<0ex,2ex>*{\cdots} \restore
\save"2,5"+<0ex,2ex>*{-2} \restore
\save"2,6"+<0ex,2ex>*{-2} \restore
\save"2,7"+<0ex,2ex>*{-2} \restore
\save"1,8"+<0ex,2ex>*{-2} \restore
\save"3,8"+<0ex,2ex>*{-2} \restore
\save"2,1"-<0ex,2ex>*{2} \restore
\save"2,2"-<0ex,2ex>*{2} \restore
\save"2,3"-<0ex,2ex>*{2} \restore
\save"2,4"-<0ex,2ex>*{\cdots} \restore
\save"2,5"-<0ex,2ex>*{2} \restore
\save"2,6"-<0ex,2ex>*{2} \restore
\save"2,7"-<0ex,2ex>*{2} \restore
\save"1,8"-<0ex,2ex>*{1} \restore
\save"3,8"-<0ex,2ex>*{1} \restore }
\\[-4mm]
& & \raisebox{-16pt}{$(ii')$} \qquad &
\xymatrix@!0@C+10pt@R-8pt{ & *=0{\bullet} \\
*=0{\circ}
  \ar@{-}[]+<1.5pt,0.9pt>;[ur]-<1.5pt,0.9pt>
  \ar@{-}[]+<1.5pt,-0.9pt>;[dr] & \\
 & *=0{\bullet}
\save"2,1"+<0ex,2ex>*{-1} \restore
\save"1,2"+<0ex,2ex>*{-2} \restore
\save"3,2"+<0ex,2ex>*{-2} \restore
\save"2,1"-<0ex,2ex>*{2} \restore
\save"1,2"-<0ex,2ex>*{1} \restore
\save"3,2"-<0ex,2ex>*{1} \restore }
\end{align*}
\vspace{-5mm}
\caption{Reducible fibres of $M'$ with self-intersections
and multiplicities.}\label{fig0}
\end{figure}
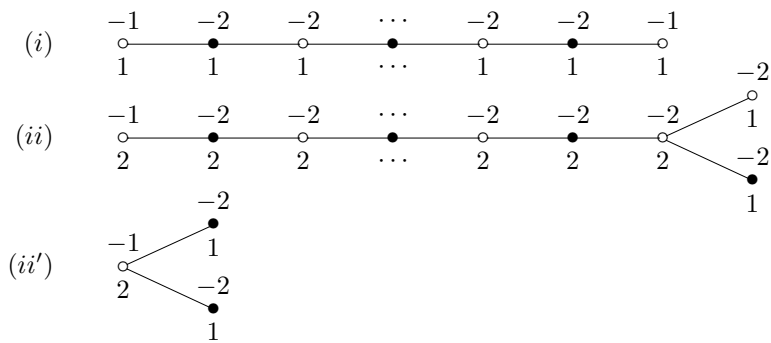

\begin{cor}\label{red}
In the above setting, if $K_{W'}^2=-2$, a reducible fibre of
$|M'|$ is of one of the following types:
\begin{enumerate} [\rm(i)]
\item either it contains two irreducible $(-1)$-curves, in which
case the dual graph of the fibre is a chain as in Figure \ref{fig0}, (i);

\item or it contains only one $(-1)$-curve with multiplicity 2, in
which case the dual graph of the fibre is shown in Figure \ref{fig0}, (ii) and (ii').
\end{enumerate}
In Figure \ref{fig0} the black (resp.\ white) vertices
correspond to curves contained (resp.\ not contained) in the branch locus.
\qed
\end{cor}

\vspace{-3mm}
\begin{rem}\label{pencil2}
Note that the number of the vertices in the graphs in Figure
\ref{fig0} is always odd.
For every $i=1, \ldots, 5$, one has $C'_i M'=C'_i (K_{W'}+D')=0$,
thus the $(-2)$-curve $C'_i$ is contained in a curve of the pencil
$|M'|$.
It is easy to verify (cf.\ also \cite{horikawa}) that
the black vertices in Figure \ref{fig0}
correspond to the curves $C'_1,\ldots,C'_5$, except
one of the $(-2)$-curves in (ii'), which is an
irreducible component of $B'_0$.

Since $M'D'=2$, the pencil $(f\circ \tpi)^* (|M'|)$ is  a genus 2
pencil without base points on  $V$ which descends via $\epsi$ to a
genus 2 pencil without base points on S.
\end{rem}

\begin{prop}\label{a=1,3}
In the above setting, if $K_{W'}^2=-2$, then there exists a
birational morphism $g:W'\to \F_a$ with either $a=1$ or $a=3$. Let
$\Delta$ be the exceptional divisor of $g$, let $E$ be the
$(-a)$-curve on $\F_a$ and $\Theta=g^*(E)$ its total transform on
$W'$. Then
\begin{equation}\label{D'}
D'\equiv 2\Theta+(a+3)M'-\Delta.
\end{equation}
Furthermore, if $a=3$, then $\Theta$ is also the proper transform
of $E$ on $W'$.
\end{prop}

\begin{proof}
As seen in Remark \ref{pencil2}, each  $(-2)$-curve $C'_i$ is
contained in a (reducible) curve of the pencil $|M'|$.

We claim that there exists a birational morphism $g:W'\to\F_a$,
for some $a\ge0$, which contracts each of the curves $C'_i$ to
points. Suppose that the fibre to which $C'_i$ belongs is of type
(i) of Corollary \ref{red}. In a birational morphism of $W'$ to a
$\F_a$, every component of such a fibre is contracted to a point,
except one of the components corresponding to an end-point.
This proves our claim in this case.

Otherwise, the fibre is of type either (ii) or (ii') of Corollary
\ref{red}. In a birational morphism of $W'$ to a $\F_a$, every
component of such a fibre is contracted to a point, except either
one of the two $(-2)$-components of multiplicity 1, corresponding
to one of the two right end-points of the graphs
in Figure \ref{fig0}.
In both cases, there is only one of the
$(-2)$-components which is one the $C'_i$'s (cf.\ Remark
\ref{pencil2}), thus we may and will choose to contract it. This
concludes the proof of our claim.

Recall that, since $K_{\F_a}\equiv -2E-(a+2)F$, where $F$ is the
ruling of $\F_a$, $K_{W'}\equiv -2\Theta-(a+2)M'+\Delta$. Then
formula \eqref{D'} follows from $ K_{W'}+ D'\equiv M'$.

Since $D'$ is nef, one has that $0 \le D' \Theta =-2a + a + 3$ by
\eqref{D'}, hence $a\le 3$.

Let $\bar\Theta$ be the proper transform of $E$ on $W'$. Again,
since $D'$ is nef, one has
\begin{equation}\label{DeltaTheta}
0\le D' \bar\Theta = -2a+a+3-\Delta\bar\Theta
 = -a+3-\Delta\bar\Theta.
\end{equation}
In particular, if $a=3$, it follows that $\Delta\bar\Theta=0$,
which means that $\bar\Theta=\Theta$.

It remains to prove that $a$ is odd. Notice that $g$ induces
a double cover $\bar\pi: \bar V \to \F_a$ branched along
the curve $\bar B_0=g(B'_0)$, which has to be an even divisor in
$\Pic(\F_a)$. Since
\begin{equation}\label{B'0a}
B'_0\equiv D'-2K_{W'}\equiv 6\Theta+(3a+7)M'-3\Delta
\end{equation}
then $\bar B_0\equiv 6E+(3a+7)F$, which shows that $a$ has to be
odd.
\end{proof}

\begin{rem}\label{3,3}
Passing from $W'$ to $\F_a$ creates only essential
singularities of the branch curve $\bar B_0$ of the double cover
$\bar\pi:\bar V\to\F_a$. More precisely, a
fibre of type (i), as in Figure \ref{fig0}, with $2l+1$ vertices
produces $2l+1$ a point $x$ of type $[3^{2l}]$, whose tangent direction
is different from the tangent direction of the fibre $F_x$ of
$\F_a$ through $x$.
Here we denote by $[3^{n}]$ the sequence $[3,\ldots,3]$ of length $n$,
so a point of type $[3^{n}]$ is a triple point
with an infinitely near point of type $[3^{n-1}]$.

Instead, a fibre as in Figures \ref{fig0},
(ii) and (ii'), with $2l+1$ vertices produces again a point $x$ of type
$[3^{2l}]$, but the tangent direction coincides with the one of
$F_x$. In case (ii') (resp.\ case (ii)),
the fibre $F_x$ is (resp.\ is not) a component of $\bar B_0$.
\end{rem}

Now we show that the case $a=3$ can be reduced to the case $a=1$.

\begin{lem}\label{l:a=3}
If $a=3$, the branch curve $\bar B_0$ of the
double cover $\bar\pi:\bar V\to\F_3$ is of
the form $\bar B_0=E+\bar B$, where $E$ is the $(-3)$-curve in
$\F_3$ and $\bar B E=1$.
Furthermore, there is a birational map $h:\F_3\map\F_1$ such that
$h\circ g:W'\to \F_1$ is a birational morphism contracting all the
$(-2)$-curves $C'_i$, $i=1,\ldots,5$, to points as in Proposition
\ref{a=1,3}.
\end{lem}

\begin{proof}
Formula \eqref{D'} says that $\bar B_0\equiv6E+16F$, which implies
the first assertion.

By Remark \ref{3,3}, $\bar B$ has five points of type $[3,3]$,
which are off $E$, because $E \bar B=1$, and can be proper or
infinitely near. Let $x$ be a $[3,3]$-point of $\bar B_0$ and let
$F_x$ the fibre of the ruling of $\F_3$ through $x$. Since $\bar B
F_x=5$, again Remark \ref{3,3} implies that the fibre $F_x$ is not
of type (ii) as in Figure \ref{fig0}. If $F_x$ is of type (ii'),
then $F_x$ is a component of $\bar B$ and
$\bar B E = F_x E = 1$, thus there is at most one such fibre.
Therefore there are fibres of type (i) and
we may choose a point $x\in\F_3$ such that $\bar B$ has a
$[3,3]$-point at $x$, whose tangent direction is different from
the one of the fibre $F_x$ through $x$ (and $F_x$ is not a
component of $\bar B$).

Now we define the birational map $h:\F_3\map\F_1$ as the
composition of two elementary transformations:
\begin{itemize}
\item the first one is based at the point $x$. This leads to
$\F_2$ and the proper transform of $\bar B$ shows the triple point
$y$ which was infinitely near to $x$; \item the second elementary
transformation is based at $y$.
\end{itemize}
The birational map $h\circ g$ is a morphism because the
exceptional curves created by the two elementary transformations
are contracted by $g$, by construction.
\end{proof}

\begin{rem}\label{r:a=3}
In the proof of the above lemma, we saw that, if $a=3$, then there
is no fibre of type (ii), and at most one fibre of type (ii').

Since $\bar B_0$ has a double point where $\bar B$ meets $E$, then
$\bar B$ has at most one unessential singularity, which is
resolved by one blowing-up (cf.\ Remark \ref{KW>-4}).
\end{rem}

Finally we deal with the case $a=1$ of Proposition \ref{a=1,3}.

\begin{lem}\label{mult5}
In the above setting, if $a=1$ and $E$ is a component of the
branch locus $\bar B_0$ of $\bar\pi:\bar V\to\F_1$, then $\bar
B_0=E+\bar B$ where $\bar B$ has a tacnode at a point $x\in E$
where the tangent direction is also tangent to $E$ at $x$.
Moreover $\bar B$ meets $E$ transversally at another point $y\ne
x$.
\end{lem}

\begin{proof}
By formula \eqref{B'0a}, one has $E \bar B=5$. Since $\bar B_0$
has at most two unessential singularities (cf.\ Remark
\ref{KW>-4}), then a $[3,3]$-point of $\bar B_0$ must lie on $E$.
\end{proof}

Let now $g':\F_1\to\pp^2$ be the contraction of the $(-1)$-curve $E$ of $\F_1$ to a
point $p\in\pp^2$.
This induces a double plane $\pi^\bullet:V^\bullet \to \pp^2$
branched along $B^\bullet_0=g'(\bar B_0)$.

The following corollary proves the Classification Theorem
in the case $K_{W'}^2=-2$.

\begin{cor}\label{Camp1}
Let $S$ be a numerical Godeaux surface with an involution
and suppose that $K_{W'}^2=-2$.
Then the double plane $\pi^\bullet:V^\bullet\to\pp^2$, constructed as above,
is branched along a reduced curve $B^\bullet_0$ of degree 10 with a point
$q$ of multiplicity at least 4 and five points of type $[3,3]$,
which can be distinct or infinitely near, not lying on a conic,
and possibly further unessential singularities resolved by at most
two blowing-ups.
\end{cor}

\begin{proof}
The curve $B^\bullet_0$ has degree 10 by \eqref{B'0a} and
multiplicity either 5 or 4 at $q$, depending on whether $E$ is or
is not a component of $\bar B'_0$. By Remark \ref{3,3},
$B^\bullet_0$ has five $[3,3]$-points, which can be proper or
infinitely near. The unessential singularities of $B^\bullet_0$
are at most two, because $K_{W'}^2-K_W^2= -2-K_W^2 \le 2$ (cf.\
Remark \ref{KW>-4}).

The fact that there is no conic through $q$ and the five
$[3,3]$-points follows by standard double plane considerations,
since $p_g(S)=0$.
\end{proof}

\begin{rem}
Lemma \ref{mult5} completely describes the nature of the
singularity at $q$ in the case $B^\bullet$ has a point of
multiplicity 5 at $q$ (cf.\ also Example 13.3 of \cite{CF}).
\end{rem}

\begin{rem} \label{r:camp}
Notice that, if $B^\bullet$ is irreducible, then its geometric
genus is 0. Examples of this type can be found in \cite{Re2},
\cite{St}, \cite{CG} (as shown by \cite{DW}).
\end{rem}

\section{On numerical Godeaux surfaces of Du Val type.}\label{duval}

In this section we go on using the notation introduced
above and we study the remaining case $K_{W'}^2=-1$
(cf.\ Lemma \ref{K2W'}).
We will show that it gives rise
to numerical Godeaux surfaces which are birationally equivalent
to double planes branched along a curve of degree 14
as in case (2) of the Classification Theorem, stated in the introduction.

According to the terminology introduced in Example 3.7, (c)
of \cite{SantaCruz},
this double plane is a
degeneration of a Du Val \emph{ancestor} with invariants $p_g=4$ and $K^2=8$.
As we said in the introduction, we will call these surfaces of \emph{Du Val type}.

We start by studying the linear systems $|M'|$  and  $|K_{W'}+M'|$.
Recall that by definition $|M'| =|K_{W'}+D'|$
and $|M'|$ has no fixed part by Lemma \ref{partfix}.

\begin{lem} \label{l:pari}
Let $S$ be a numerical Godeaux surface with an involution.
and let $W'$ as in Sections \ref{prima}, \ref{sGodeaux}.
Suppose that $K_{W'}^2=-1$. Then:
\begin{enumerate}[\rm (i)]
\item $|M'|$ is a pencil of elliptic curves with a simple
base point;
\item $h^0(W',\OO_{W'}(K_{W'}+M'))=1$;
\item if $G$ is the unique curve in $|K_{W'}+M'|$,
then $G M'=0$ and $G=2E+C'_5$, where
$E$ is a $(-1)$-curve and $E C'_5=E D'=1$;
\item the divisor $C'_1+\cdots+C'_4$ is even in $\Pic(W')$.
\end{enumerate}
\end{lem}

\begin{proof}
Since $K_{W'}^2=-1$, one has $M'^2=(K_W'+D')^2= 1$. Moreover
$K_{W'}M'= -1$. This proves (i).
The long exact sequence obtained  from
\[
0 \to \OO_{W'}(K_{W'}) \to \OO_{W'}(K_{W'}+M') \to \OO_{M'} \to 0
\]
implies part (ii), because $W'$ is rational.

Let us prove part (iii).
Since $p_a(M')=0$, one has that $M'G=0$
and every component $\theta$ of $G$ is such that $M'\theta=0$, because $M'$ is nef.
Furthermore  $M' G=0$ implies that  the
intersection form on the components of  $G$ is negative definite. Since $G^2=-2$,  there is
an irreducible curve $E$ in $G$ such that $E^2<0$ and
$E(K_{W'}+M')=E G<0$. From $M'E=0$, we conclude that   $E K_{W'}<0$, thus $E$ is a $(-1)$--curve and $E D'=1, EG=-1$.
Furthermore $E(G-E)=0$.

We claim that $E$ intersects one of the $(-2)$--curves  $C'_i$, say $C'_5$.  Indeed, since
$D'+C'_1+\cdots+C'_5 $ is divisible by $2$ by Remark \ref{divisibilita},
the claim follows, because $E D'=1$.

Note that
$C'_5 G=C'_5 (2K_{W'}+D')=0$. Since
$E C'_5>0$, we have $C'_5\leq G-E$.
Since $E(G-E)=0$, then $E$ is necessarily a component of $G-E$.
Note that $(2E+C'_5)^2=4E^2+4EC'_5+C'^2_5=4EC'_5-6  $ and so,
because the intersection form on the components of $G$ is negative definite,
we conclude that  $(2E+C'_5)^2=-2$.
Since $ (2E+C'_5) G=-2$, we have $(G-(2E+C'_5))^2=0$ and,
so again  by negative definiteness,  we see that $G= 2E+C'_5$.

So the proof of part (iii) is concluded.

Finally we prove part (iv). Since $D'+C'_1+\cdots+C'_5$ is
divisible by $2$, then
$2K_{W'}+D'+C'_1+\cdots+C'_5 \equiv C'_5+2E+C'_1+\cdots+C'_5$ is
divisible by $2$ and therefore also $C'_1+\cdots+C'_4$ is divisible by $2$.
\end{proof}

\begin{rem}\label{r:pari}
Notice that $C_1+\ldots+C_4$
is an even divisor in $W$ by the previous lemma,
hence $B_0+C_5$ is even too.
So, by Beauville's Lemma \ref{l:2-tors},
$\Tors_2(S)$ is not trivial
and thus, by Theorem \ref{t:tors}, $\Tors(S)$ is either $\Z_2$ or $\Z_4$.
\end{rem}

We will need later the following:

\begin{cor} \label{connected}
In the above setting,
the curve $B_0'$ is 1-connected.
\end{cor}

\begin{proof}
Suppose that $B_0'$ is not $1$-connected
and write $B_0'=G+H$ with $GH=0$.
Then say $D'G=0$ and $D'H=2$.
Since $D'$ is nef, $G^2<0$  and
therefore $0=D'G=2K_{W'}G+GH+G^2= 2K_{W'}G+   G^2$
implies that $K_{W'}G>0$.
Also we note that $H^2\leq 1$.
In fact the Index Theorem applied to $D'$ and $D'-H$
implies that $H^2\leq 2$ and that $H^2=2$ if
and only if $D'\equiv H$ (because $W'$ is rational).
But if $D'\equiv H$, then $G\equiv -2K_{W'}$
and, since also $C_1+\cdots+C_4$ is even,
this would imply by Beauville's Lemma \ref{l:2-tors} that  $|\Tors_2(S)|\geq 2$,
contradicting Theorem \ref{t:tors}.
Then $2=D'H=  2K_{W'}H+   H^2$ implies that also  $K_{W'}H>0$.
Since $K_{W'}B_0' =2$, it follows that $K_{W'}H=K_{W'}G=1$.
From $0=D'G=2K_{W'}G +   G^2$, we conclude that $G^2=-2$,
which contradicts the adjunction formula.
\end{proof}

By Lemma \ref{l:pari}, one has that
\begin{equation} \label{KW'}
K_{W'}\equiv -M' + 2E+C'_5.
\end{equation}
Furthermore, since $C'_1+\cdots+C'_4$ is an even divisor,
also $D'-C'_5\equiv 2\Delta$, where $\Delta^2=0$ and $K_{W'}\Delta=0$.
By the Riemann-Roch Theorem, $|\Delta|\neq \emptyset$.
The following identities are easy to check:
\begin{align}
& M'\equiv \Delta+E+C'_5, \label{M'} \\
& B'_0\equiv 4\Delta-2E+C'_5\equiv 4M'-6E-3C'_5, \label{B'0equiv} \\
& K_{W'}+\Delta\equiv E. \label{KW'+Delta}
\end{align}

\begin{rem}\label{r:tildeDelta}
Going back to our original surface $W$, the above formulas
mean that $D\equiv 2\tilde \Delta+ C_5$,
where $\tilde\Delta$ is the pull-back of $\Delta$ to $W$.
\end{rem}

Note that there is a birational morphism $g:W'\to X$ which
contracts $G=2E+C'_5$ to a smooth point $q$.
The rational surface $X$ is such that $K_X^2=1$ and $g^*(K_X)\equiv -M'$.
Therefore $-K_X$ is nef and big.

We denote by $\DD$ the pencil $g_*(|D'|)$ and by  $D''\in\DD$ its
general element. Then $D''^2=4$, $K_XD''=-2$ and $\DD$ has at $q$ a base point
with a fixed tangent. Furthermore  $D''\equiv -2K_X$.

We note that $X$ still contains an even set of four  disjoint $(-2)$-curves.
Now we can apply the following result from \cite{CCM}, regarding
rational surfaces with an even set of nodes
(cf.\ also \S4 in Chapter 0 of \cite{CD}).

\begin{thm}\label{t:delPezzo}
Let $X$ be a weak Del Pezzo surface, i.e.\ $-K_X$ is big and nef.
Assume that $K_X^2=1$ and that $X$ has four disjoint $(-2)$-curves
$C''_1,\ldots,C''_4$ such that $C''_1+\cdots+C''_4$ is even
in $\Pic(X)$.
Then there exists a birational morphism $h:X\to X'$,
where $X'$ is obtained from $\F_a$, with $a=0$, $1$, or $2$, by blowing up:
\begin{itemize}
\item two points $q'_1$, $q'_2$ in distinct fibres $F_1,F_2$ of the same ruling of $\F_a$;
\item the point which is the intersection
of the strict transform of $F_i$ with the exceptional curve corresponding to $q_i$, $i=1,2$;
\end{itemize}
in case $a=2$, none of the blown-up points lies on the $(-2)$-curve of $\F_2$.
The morphism $h$ maps isomorphically the $(-2)$-curves $C''_1,\ldots,C''_4$
onto the proper transforms in $X'$ of $F_1$, $F_2$
and of the exceptional curves corresponding to $q'_1$ and $q'_2$.
\qed
\end{thm}

Next corollary proves the Classification Theorem in the last case (2).

\begin{cor}\label{c:duval}
Let $S$ be a numerical Godeaux surface with an involution
and suppose that $K_{W'}^2=-1$.
Then there exists a birational morphism $g':W'\to\pp^2$ such that
the induced double cover $\pi^\bullet:V^\bullet\to \pp^2$
is a double plane of Du Val type,
branched along a reduced curve $B^\bullet$, which possibly
has irrelevant singularities, resolved with at most three blowing-ups.
\end{cor}

\begin{proof}
Suppose that $a=1$ in the statement of Theorem \ref{t:delPezzo},
i.e.\ there exists a birational morphism $X\to\F_1$.
By blowing-down the $(-1)$-curve of $\F_1$, it determines
a birational morphism $h':X\to\pp^2$ and hence
a birational morphism $g'=h'\circ g\circ f:W\to\pp^2$.
Up to reordering the indices, the $(-2)$-curves $C_1$ and $C_2$
are mapped via $g'$ to two distinct lines, say $r_1$ and $r_2$,
whilst $C_3$ and $C_4$ are contracted to two distinct points
$q_1\in r_1$ and $q_2\in r_2$.
The anti-canonical pencil $|-K_X|$ on $X$ is mapped via $h'$ to a pencil
of plane cubics with the following eight base points:
\begin{itemize}
\item $q_0=r_1\cap r_2$;
\item $q_i$, $i=1,2$, where the cubics are tangent to $r_i$ at $q_i$;
\item further three simple base points $q_3,q_4,q_5$.
\end{itemize}
Therefore the pencil $|D|$ on $W$ is mapped via $g'$
to a pencil of plane sextics with the following base points:
\begin{itemize}
\item double points at $q_0,q_3,q_4,q_5$;
\item a tacnode at $q_i$, $i=1,2$, where the tacnodal tangent is $r_i$;
\item a further simple base point $q_6$ with fixed tangent.
\end{itemize}
Standard double plane considerations imply the assertion in this case $a=1$.

Suppose that $a=2$ in the statement of Theorem \ref{t:delPezzo},
i.e.\ there exist birational morphisms $X\to X'\to\F_2$.
Note that on $X'$ there is a $(-2)$-curve $N$, which is the proper transform of
the $(-2)$-curve on $\F_2$, and two reducible fibres of  type $A_1+A_3+2N_1$
and $A_2+A_4+N_2$,
where $A_1,\ldots,A_4$ are $(-2)$-curves and $N_1,N_2$ are $(-1)$-curves with
$N A_1=N A_2=1$, $N A_3= N A_4 = N N_1 = N N_2 = 0$.
By first blowing-down $N_1$, $N_2$ and then the image of $A_1$ and of $A_4$,
one arrives to $\F_1$ and the image of $N$ is the $(-1)$-curve.
At this point, one proceeds as in the previous case.
Notice that, in the present situation, either $q_1$ or $q_2$ is infinitely
near to $q_0$.

Suppose finally that $a=0$.
On $X'$ there are two reducible fibres of the same ruling of the type
$A_1+A_3+2N_1$ and $A_2+A_4+2N_2$, where $A_1,\ldots,A_4$ are $(-2)$-curves
and $N_1,N_2$ are $(-1)$-curves such that, if $N_0$ is a general fibre
of the other ruling,
$N_0 A_1=N_0 A_2=1$, $N_0 A_3= N_0 A_4 = N_0 N_1 = N_0 N_2 = 0$.
By first blowing-down $N_1$, $N_2$ and then the image of $A_1$ and of $A_4$,
one arrives to $\F_1$, because the image of $N_0$ is a section with self-intersection 1,
and one concludes as before.
\end{proof}

\begin{rem} \label{r:duVal}
Notice that, if the degree 12 component of $B^\bullet$ is irreducible,
then its geometric genus is 1.

Remark also that all irrelevant singularities are double points,
because four blowing-ups are needed to resolve an irrelevant triple point
(cf.\ Remark \ref{KW>-4}).
\end{rem}

\begin{rem} \label{r:Delta}
The morphism $g':W\to\pp^2$ defined in the proof of Corollary \ref{c:duval}
maps the curve $\tilde\Delta$ of Remark \ref{r:tildeDelta}
to a plane cubic $\bar\Delta$.
Indeed the curve $\Delta$ on $W'$, cf.\ formulas \eqref{KW'}--\eqref{KW'+Delta},
is mapped to $\bar\Delta$ via $h'\circ g:W'\to\pp^2$.
\end{rem}

\begin{rem} \label{Gamma}
Consider $X'$ as in Theorem \ref{t:delPezzo}.
Let $\bar C_1$, $\ldots$, $\bar C_4$ be the $(-2)$-curves of $X'$
such that $\bar C_1+\cdots+\bar C_4$ is even.
Let $\bar\Upsilon$ be the proper transform of a general conic through $q_1$ and $q_2$
and suppose that $q_i$, $i=1,2$, is not infinitely near to $q_0$.
Clearly $-2K_{X'}\equiv 2\bar\Upsilon + \bar C_1+\cdots+\bar C_4$.
Going back to $W$, setting
$H$ and $\Upsilon$ the pull-back of $-K_{X'}$ and $\bar\Upsilon$, respectively,
then $C_1+\cdots+C_4\equiv 2(H-\Upsilon)$.
\end{rem}

In Corollary \ref{c:godeaux} we saw that if $R$ has
an irreducible component $\Gamma$ of genus 2,
then $W$ is either an Enriques surface or a rational surface
with $K_W^2=-2$ or $-1$.
Now we conclude this section by ruling out the case that $W$ is rational,
as announced.

\begin{prop} \label{p:exc}
Let $S$ be a numerical Godeaux surface with an involution
(cf.\ notation in \S\ref{sGodeaux}).
Suppose that $-2 \leq K_W^2 \leq -1$ and $\Gamma^2=1$.
Then $W$ is not rational.
\end{prop}
\begin{proof}
Suppose by contradiction that $W$ is rational. Let $\Gamma_0=\eta^*(\pi\Gamma)$.
By the proof of Corollary \ref{c:godeaux},
one has that $B_0-\Gamma_0\equiv -2K_W$.
Beauville's Lemma \ref{l:2-tors} and Theorem \ref{t:tors} imply that
in this case $\sum_{i=1}^4 C_i$ cannot be divisible by 2.

By the above Lemma \ref{l:pari}, it follows that $W'=W$ and $K_{W}^2=-2$.
In this case $B_0=\Gamma_0+N_1+N_2$, where $N_1,N_2$ are smooth rational curves
with self-intersection $-4$ by formula \eqref{e:h}.
Then, because  $B_0-\Gamma_0\equiv -2K_W$,  $N_1+N_2\equiv -2K_W$.
Note that $D-N_1-N_2\equiv D+2K_W$ cannot be effective because
$|M|=|D+K_W|$ is a base point free pencil of rational curves by Lemma \ref{l:M'rational}.

Since $\Gamma_0\in|D|$ is irreducible,
the pencil $|D|$ has no fixed part.
In particular any two distinct curves of $D$ have no common component.
Since $D N_1=D N_2=0$, we can write $D\equiv G+N_1 \equiv H+N_2$.
Since, as we saw, $D-N_1-N_2$ is not effective,
the two curves $G+N_1$ and $H+N_2$ in $|D|$ are distinct,
hence $G,H$ are effective divisors without common components. Note that $G^2=H^2=-2 $ and $K_WG=K_WH=-2$
Since $2D\equiv G+H+N_1+N_2\equiv G+H-2K_W$,
one has $G+H\equiv 2(D+K_W)\equiv 2M$.

By Lemma \ref{l:M'rational}, every component $A$ of $G$ is either a $(-1)$-curve
or a $(-2)$-curve.
In the latter case, $0=A D=A H + A N_2$ implies that
$A H=A N_2=0$, because $G$ has no common components with $H$. Since $N_2\equiv -2K_W-N_1$, also $A N_1=0$, implying that $A G=0$.
Moreover, since $K_W G=-2$, one has that $G$ contains either two $(-1)$-curves
$A_1,A_2$ or one $(-1)$-curve $A$ with multiplicity 2.
In the first case, $A_i G\geq A_i^2=-1$,
hence, since $G^2=-2$, the equality holds
and $(G-A_1-A_2)^2=0$, implying by Zariski's Lemma that $G=A_1+A_2$.
In the second case, one sees that $A G=-1$, implying that
$A$ meets $G-2A$ transversally at one point of a $(-2)$-curve $\theta$.
As above we conclude that $G=2A+\theta$. Similarly one shows that either $H=A'_1+A'_2$ or $H=2A'+\theta'$,
where $A',A'_i$ are $(-1)$-curves and $\theta'$ is a $(-2)$-curve.

Now recall that $2M\equiv G+H$, $G$ and $H$ have no common components and $|D-M|=\emptyset$ . So $G+H=M_1+M_2$ where $M_1, M_2$ are distinct fibres of $|M|$ and this   excludes the possibility that $G=2A+\theta$
or $H=2A'+\theta'$.

Therefore the only remaining possibility is that $G=A_1+A_2$ and $H=A'_1+A'_2$.
This implies that $M_1=A_1+A'_1$ and $M_2=A_2+A'_2$,
after possibly reordering the indices, which is impossible by Corollary \ref{red}.
\end{proof}

\section{On the torsion of numerical Godeaux surfaces of Du Val type.}\label{s:torsion}

In this section we will study the torsion of a numerical Godeaux surface $S$
of Du Val type.
We will freely use the notation of the previous sections.

\begin{prop} \label{p:tors}
Let $S$ be a numerical Godeaux surface of Du Val type.
Then there is a non-trivial 2-torsion element in $\Tors(S)$,
and accordingly, there is an {\'e}tale 2-to-1 cover $\bar S\to S$, where
$\bar S$ is a regular surface with $p_g(\bar S)=1$ and $K^2_{\bar S}=2$.
Therefore $\Tors(S)$ is either $\Z/2\Z$ or $\Z/4\Z$,
depending on whether $\bar S$ has either no torsion
or 2-torsion.
\end{prop}

\begin{proof}
Recall that, by Remark \ref{r:pari},
the $(-2)$-curves $C_1,\ldots,C_4$ in $W$ are
such that $C_1+\cdots+C_4$ is even and so $S$ has non-trivial 2-torsion,
which defines an {\'e}tale cover $\tilde S\to S$.
The rest of the assertion follows from Theorem \ref{t:tors}.
\end{proof}

We can now consider the following commutative diagram:
\[
\xymatrix{%
 \tilde V \ar[d]^{\Pi} \ar[r]
 & V \ar[d]^{\tilde\pi} \\
\tilde W \ar[r]^{t}  & W
}
\]
where $t:\tilde W\to W$ is the double cover branched over $C_1+\cdots+C_4$,
$\Pi:\tilde V\to \tilde W$ is the double cover
branched over $t^*(B_0+C_5)$ and $\tilde V\to V$
is the \'etale double cover associated to the 2-torsion of $V$.

Note that the minimal model of $\tilde V$ is obviously isomorphic
to the surface $\tilde S$ of Proposition \ref{p:tors}.
Also remark that   any smooth rational curve $\theta$ on $W$
disjoint from the  curves $C_1,\ldots,C_4$ pulls back on $\tilde W$
to two disjoint rational curves with the same self-intersection number as $\theta$.

Note also that $t^*(B_0+C_5)= C_{5,1}+C_{5,2}+\tilde B_0$,
 where $\tilde B_0$ is an {\'e}tale double cover
of $B_0$, which is  disjoint from the two $(-2)$-curves $C_{5,1}$ and $C_{5,2}$.
Furthermore, since $B_0+C_5$ is divisible by $2$ in $\Pic(W)$,
also $ C_{5,1}+C_{5,2}+\tilde B_0 \equiv 2\Phi$ with $\Phi$ in $\Pic (\tilde W).$

\begin{lem} \label{B''_0reducible}
Let $S$ be a numerical Godeaux surface of Du Val type, and
suppose that $\Tors(S)=\Z/4\Z$.
Then the double cover $t$ splits over $B_0+C_5$.
\end{lem}

\begin{proof}
By Proposition \ref{p:tors},
$\tilde V$ has 2-torsion if and only if $\Tors(S)=\Z/4\Z$.
On the other hand, $\tilde W$ is rational and thus,
by Beauville's Lemma \ref{l:2-tors}, $\tilde V$ has 2-torsion
if and only if there exists a curve
$G\lneqq\tilde B_0+C_{5,1}+C_{5,2}$
divisible by $2$ in $\Pic(\tilde W)$.
Let $\iota$ be the involution on $\Pic(\tilde W)$ (and $H^2(\tilde W, \Z)$)
induced by the involution on $\tilde W$ corresponding to the double cover $t$.
We note that $\Phi$, where $ C_{5,1}+C_{5,2}+\tilde B_0\equiv 2\Phi$,
is invariant under $\iota$.

Suppose that  there exists a curve
$G\lneqq\tilde B_0+C_{5,1}+C_{5,2}$ divisible by $2$ in $\Pic(\tilde W)$.
Then also $ H:=\tilde B_0+C_{5,1}+C_{5,2}-G$ is divisible by $2$
and we can write $G\equiv 2L_1$, $H\equiv 2L_2$ where $L_1+L_2\equiv\Phi$.
Note that neither $L_1$ nor $L_2$ can be invariant under $\iota$,
because otherwise there would be too many divisibility relations
in the branch locus in $W$, implying again by Beauville's Lemma \ref{l:2-tors}
that $|\Tors_2(S)|\geq 2$,
which is impossible by Theorem \ref{t:tors}.

Write then $G=G'+J_1$, and $H=H'+J_2$
where all divisors appearing are effective
and we assume that $\iota(G')=G'$, $\iota (H')=H'$, $\iota(J_i)\neq J_i$, $i=1,2$
(and $\iota(\theta)\neq \theta$ for any component $\theta$ of $J_i$).
Since $\iota(G+H)=G+H$, necessarily $\iota(J_1)=J_2$.
Note also that $J_i\neq 0$, $i=1,2$,
otherwise $\iota(G)=G$ and so $ \iota(2L_1)=2\iota (L_1)=2L_1$,
implying, because $\tilde W$ is rational, that $\iota(L_1)=L_1$.
This would mean that $\iota (G)=G$ and the divisibility would be coming
already from $W$, which is impossible.

Now we have $G+\iota(G)= G'+J_1+G'+J_2=2G'+J_1+J_2$.
On the other hand, $G+\iota(G) \equiv 2L_1+2\iota(L_1)$
and so we conclude that $J_1+J_2 \equiv 2L_1+2\iota(L_1)-2G'$ is divisible by $2$.
Again by Beauville's Lemma \ref{l:2-tors} and Theorem \ref{t:tors},
there cannot exist further divisibility relations in $\tilde B_0+C_{5,1}+C_{5,2}$,
hence we conclude that $G'=H'=0$, $\iota(L_1)=L_2$, $\iota (J_1)=J_2$
and therefore $t$ splits over $B_0+C_5$.
\end{proof}

\begin{lem} \label{split}
Let $S$ be a numerical Godeaux surface of Du Val type and
suppose that  the double cover $t$ splits over $f^*(B'_0+C'_5)$,
where $f:W\to W'$ is the birational morphism of Proposition \ref{W'}.
Then $\Tors(S)=\Z/4\Z$.
\end{lem}

\begin{proof}
Note that the components of $f^*(B'_0)$, which are not components of $B_0$,
form trees of rational curves; each one of them is double in $f^*(B'_0)$
and does not meet $C_1,\ldots,C_4$, so $t$ splits over it.

Write $t^* (f^*(B'_0+ C_5')) = A_1 +C_{5,1}+A_2+ C_{5,2}$
where $\iota (A_1)=A_2$ and $A_i$, $i=1,2$, is isomorphic to $f^*(B'_0)$,
so it is connected by Corollary \ref{connected}
and each component of $A_i$, which is not a component of $t^*(B_0)$,
is double in $A_i$.
Recall that $B_0'+C_5'\equiv 4\Delta-2E+2C'_5$ by formula \eqref{B'0equiv}.
Denote the two connected components of $t^* (f^*(E))$
by $E_1$ and $E_2$, where $E_1C_{5,1}=1$ and  $E_1C_{5,2}=0$.
Since $E B_0'=3$ and $(E+B_0')^2>0$,  $t^* (f^*(E+B_0'))$ is connected
and therefore, say, $E_2A_1=2$ and $E_2A_2=1$.

Then we claim that $A_1\equiv  t^* (f^* (2\Delta))-2E_2+C_{5,1}$.
In fact we have $t^* (D) A_1= 2$  and  $t^* (D)(t^* (f^* (2\Delta))-2E_2+C_{5,1})=2$.
Since $(A_1- t^* (f^* (2\Delta))+2E_2-C_{5,1})^2=0$,
by the Index Theorem we obtain the claim, because $\tilde W$ is a rational surface.

Hence $A_1+C_{5,1}$ is divisible by $2$ in $\Pic(\tilde W)$ and therefore also
$t^*(B_0+C_5)$ strictly contains an effective divisor divisible by 2,
implying by Beauville's lemma \ref{l:2-tors}
that $\tilde V$ (and $\tilde S$) has 2-torsion.
\end{proof}

\begin{rem}
Note that if $t$ splits over $f^*(B'_0+C'_5)$, then $t$ splits also over $B_0+C_5$,
but the converse is not necessarily true.
\end{rem}

For the examples of the next section,
we will need the following criterion:

\begin{thm} \label{settica}
Let $S$ be a numerical Godeaux surface $S$ of Du Val type
(cf., e.g., Cor\-ol\-lary \ref{c:duval} and notation therein).
Suppose that $B_0$ on $W$ is a (smooth) irreducible curve with $p_a(B_0)=1$,
and that $q_i$, $i=1,2$, is not infinitely near to $q_0$.

Then $\Tors(S)=\Z/4\Z$ if and only if
there exists a (unique) plane curve of degree 8 with the following singularities:
\begin{itemize}
\item a point of multiplicity $2$ at $q_0$;
\item a point of type $[3,2]$ at $q_i$, $i=1,2$, where the infinitely
near double point is in the direction of the line $r_i$;
\item triple points at $q_3$, $q_4$, $q_5$;
\item a tacnode at $q_6$, where the tacnodal tangent is tangent also to $B^\bullet$.
\end{itemize}
\end{thm}

\begin{proof}
The hypotheses imply that $W'=W$.
By Lemmas \ref{B''_0reducible}-\ref{split},
$\Tors(S)=\Z/4\Z$ if and only if
the double cover $t$ splits on $B_0$, which happens by Remark \ref{Gamma} if and only if
the restriction of $-(C_1+\cdots+C_4)/2\equiv -H+\Upsilon$ to $B_0$ is trivial,
or, equivalently, if and only if
$
\OO_{B_0}(K_{W}+B_0-H+\Upsilon) \cong \OO_{B_0}.
$
Notice that, since $B_0 \cap C_i =\emptyset$, $i=1,\ldots,4$,
$\OO_{B_0}(K_{W}+B_0-H+\Upsilon)$ has degree 0.
The sequence
\[
0\! \to\! \OO_{W}(K_{W}-H+\Upsilon) \!\to\! \OO_{W}(K_{W}+B_0-H+\Upsilon)
\!\to\! \OO_{B_0}(K_{W}+B_0-H+\Upsilon) \!\to\! 0
\]
is exact.
Note that
$H^0(W,\OO_{W}(K_{W}-H+\Upsilon))=0$, because
$H^0(W,\OO_{W}(2K_{W}-C_1-\cdots-C_4))=0$,
and, since $H^2(W,\OO_{W}(K_{W}-H+\Upsilon))=0$,
the Riemann-Roch Theorem implies that $H^1(W,\OO_{W}(K_{W}-H+\Upsilon))=0$.

Thus $\OO_{B_0}(K_{W}+B_0-H+\Upsilon) \cong \OO_{B_0}$
if and only if $h^0(W,\OO_{W}(K_{W}+B_0-H+\Upsilon)) = 1$.
The morphism $g':W\to\pp^2$ (cf.\ the proof of Theorem \ref{c:duval})
maps the unique curve in $|K_{W}+B_0-H+\Upsilon|$ to a plane curve of degree 8
as in the statement.
\end{proof}

Theorem \ref{settica} can be extended, under suitable assumptions,
e.g.\ as follows.

\begin{cor} \label{c:settica}
Let $S$ be a numerical Godeaux surface $S$ of Du Val type.
Suppose that $q_i$, $i=1,2$, is not infinitely near to $q_0$
and the curve $B_0$ on $W$ has two irreducible
components, $\Gamma_0$ of genus 1 and $Z$ of genus 0, such that
$B'_0$ on $W'$ is the union of two curves meeting transversally at a point.
Then $\Tors(S)=\Z/4\Z$ if and only if
there exists a plane curve of degree 8, as in the statement of Theorem \ref{settica},
containing the image of $Z$ in the plane.
\end{cor}

\begin{proof}
Lemmas \ref{B''_0reducible} and \ref{split} again imply
that the torsion group of $S$ is $\Z/4\Z$ if and only if the double cover
$t:\tilde W\to W$ splits over $\Gamma_0$.
The same argument of the proof of Theorem \ref{settica} shows that
there is the plane curve of degree 8 as above
and it contains the rational curve
which is the image of $Z$ in the plane.
\end{proof}

\begin{rem}\label{exkeum}
By Propositions \ref{godeaux} and \ref{p:exc}
(cf.\ Remarks \ref{r:camp} and \ref{r:duVal}),
it follows that numerical Godeaux surfaces
described in Example 4.2 in \cite{keum},
which have $\Tors=\Z/4\Z$,
are birational to double planes of Du Val type.
In that example, indeed, the ramification curve $R$ is irreducible
of genus 1 and this may happen, by our classification results,
only if $W$ is rational and $S$ is of Du Val type.
\end{rem}

A natural question is whether a numerical Godeaux surface $S$
of Du Val type can be also of Campedelli type,
of course for different involutions on $S$.

In order to give a result about this problem,
we first prove the following:

\begin{prop} \label{k+etareducible}
Let $S$ be a numerical Godeaux surface of Campedelli type.
If there is a non-trivial element $\eta_2$ in $\Tors_2(S)$,
then the unique curve $\Delta_2$ in $|K_S+\eta_2|$ is reducible
and has a common component $A$ with the ramification curve $R$.
Furthermore, $A$ is rational.
\end{prop}

\begin{proof}
Note that $\Delta_2$ is fixed by the involution
and it is not contained in $R$,
because $\Delta_2$ has arithmetic genus 2 (cf.\ Corollary \ref{c:godeaux}
and Proposition \ref{p:exc}).
Suppose that $\Delta_2$ and $R$ have no common component.
Looking at the Campedelli double plane representation (cf.\ Corollary \ref{Camp1}),
the image of $\Delta_2$ is a plane curve $\Delta'$,
which has no common component with the branch curve $B^\bullet$.
Since $2\Delta_2\in |2K_S|$,
the curve $2\Delta'$ is contained in the linear system
of quartics with a double point at $p$ and passing simply through the five points
of type $[3,3]$ with tangent lines also tangent to the $[3,3]$-point.
Then $\Delta'$ should be a conic passing through all the given points,
a contradiction.
We conclude that $\Delta_2$ and $R$ has a common component $A$.
Reasoning as above, one sees that
the image of $A$ in the plane is a conic, hence $A$ is rational.
\end{proof}

\begin{cor} \label{c:notCamp}
Let $S$ be a numerical Godeaux surface of Du Val type
and $\pi^\bullet:V^\bullet\to \pp^2$ the corresponding double plane,
branched along the curve $B^\bullet$
as in Corollary \ref{c:duval}.
Let $\bar\Delta$ be the unique plane cubic curve
passing through $q_0$, $q_i$, $i=1,2$, where
it is tangent to $r_i$, $q_3$, $q_4$, $q_5$ and $q_6$.
Suppose that $\bar\Delta$ is irreducible
and does not pass through any irrelevant singularity
of $B^\bullet$.
Then $S$ has no involution which realizes it as double plane of Campedelli type.
\end{cor}

\begin{proof}
By Lemma \ref{l:pari} and Remark \ref{r:pari},
$S$ has an element $\eta_2$ of $2$-torsion.
According to Remark \ref{r:Delta},
the unique curve $\Delta_2$ in $|K_S+\eta_2|$ is mapped
to the plane cubic $\bar\Delta$.
The hypotheses on $\bar\Delta$ imply that $\Delta_2$
is irreducible, hence there is no involution on $S$
which realizes it as double plane of Campedelli type
by Proposition \ref{k+etareducible}.
\end{proof}

\section{Examples of numerical Godeaux surfaces of Du Val type.} \label{s:ex}

As we already remarked,
there was no previously known construction of
numerical Godeaux surfaces as double planes of Du Val type.

In this section we produce such a construction.
In order to do that,
one has to find a reduced curve $B^\bullet$ of degree 12
with singularities at the points $q_0, \ldots, q_6$
as described in the Classification Theorem stated in the introduction.

If one chooses the points in general position,
then one expects no curve like $B^\bullet$,
because the virtual dimension of the linear system
of curve of degree 12 with those singularities is $-2$.
However it is possible to find such irreducible curves.
One way is to look for a curve which is invariant under a linear transformation
of the plane of order 2.
This is an idea originally used by Stagnaro in order to construct
numerical Godeaux surfaces
as double planes of Campedelli type, cf.\ \cite{St} and \cite{We3}.

\begin{example} \label{ex:Z4}
Let $[x,y,z]$ be homogeneous coordinates in $\pp^2$.
Let $r_1$ be the line $x=y$
and choose the following points:
\[
q_0 =[0,0,1], \quad
q_1 =[1,1,1], \quad
q_3 =[1,0,0], \quad
q_4 =[0,1,1], \quad
q_6 =[-2,0,1],
\]
and let $r_3$ be the line $x+2z=0$, which passes through $q_6$.
The linear involution
\[
\phi: [x,y,z]\in\pp^2 \to [x,-y,z] \in \pp^2,
\]
fixes $[0,1,0]$ and all the points of the line $y=0$.
Consider the line $r_2=\phi(r_1)$ and the points $q_2=\phi(q_1)=[1,-1,1]$
and $q_5=\phi(q_4)=[0,-1,1]$.

Note that there is no conic passing through $q_1,q_2,\ldots,q_6$.

One sees that there are 49 monomials of degree 12 in $x,y,z$ which are
invariant under $\phi$, i.e.\ their degree in $y$ is even.
Hence the curves of degree 12, which are invariant under $\phi$,
form a linear system of dimension 48.
Now we impose:
\begin{itemize}
\item quadruple points at $q_0$, $q_3$ and $q_4$;
\item a point of type $[4,4]$ at $q_1$,
where the tangent direction is the line $r_1$;
\item a point of type $[3,3]$ at $q_6$,
where the tangent direction is the line $r_3$.
\end{itemize}
Note that, in this way, we are imposing also another quadruple point at $q_5=\phi(q_4)$
and another point of type $[4,4]$ at $q_2=\phi(q_1)$,
where the tangent direction is the line $r_2=\phi(r_1)$,
so the resulting curve should have the required singularities.

The interesting observation is that $q_0$ and $q_3$, instead of imposing
10 independent conditions each, impose only 6 conditions each,
because of the symmetry of the configuration.
For the same reason, the singularities at $q_6$ does not impose 12 independent
conditions as expected, but only 6.
In conclusion, the number of independent conditions we are imposing
are no more than 48, so that there is surely a curve with at least the required singularities.

Using a computer algebra program (we used Maple for this),
it turns out that there is one curve of degree 12
satisfying the above conditions and having exactly the required
singularities, and not worse than those.
Its equation is:

\vspace{-2mm}
{\scriptsize
\begin{align*}
& 15625\,{x}^{8}{y}^{4}-117650\,{x}^{8}{y}^{2}{z}^{2}+81289\,{x}^{8}{z}^{4}
+604400\,{x}^{7}{y}^{4}z-908128\,{x}^{7}{y}^{2}{z}^{3}+386672\,{x}^{7}{z}^{5}+ \\[-0.75mm]
& -867475\,{x}^{6}{y}^{6}
+1931246\,{x}^{6}{y}^{4}{z}^{2}-1557283\,{x}^{6}{y}^{2}{z}^{4}
+ 369096\,{x}^{6}{z}^{6}-561632\,{x}^{5}{y}^{6}z+ \\[-0.75mm]
& +429504\,{x}^{5}{y}^{4}{z}^{3}
+777504\,{x}^{5}{y}^{2}{z}^{5} -562432\,{x}^{5}{z}^{7}+857979\,{x}^{4}{y}^{8}
-520134\,{x}^{4}{y}^{6}{z}^{2}+ \\[-0.75mm]
& -2103701\,{x}^{4}{y}^{4}{z}^{4}+2553616\,{x}^{4}{y}^{2}{z}^{6}
-808496\,{x}^{4}{z}^{8}-572400\,{x}^{3}{y}^{8}z+2171200\,{x}^{3}{y}^{6}{z}^{3}+ \\[-0.75mm]
& -2711600\,{x}^{3}{y}^{4}{z}^{5}+1112800\,{x}^{3}{y}^{2}{z}^{7}
+802575\,{x}^{2}{y}^{10}-4569750\,{x}^{2}{y}^{8}{z}^{2}
+8348975\,{x}^{2}{y}^{6}{z}^{4}+ \\[-0.75mm]
& -6199000\,{x}^{2}{y}^{4}{z}^{6}
+1617200\,{x}^{2}{y}^{2}{z}^{8}+540000\,x{y}^{10}z
-1620000\,x{y}^{8}{z}^{3}+1620000\,x{y}^{6}{z}^{5}+ \\[-0.75mm]
&-540000\,x{y}^{4}{z}^{7}-810000\,{y}^{12}+3240000\,{y}^{10}{z}^{2}-4860000\,{y}^{8}{z}^{4}
+3240000\,{y}^{6}{z}^{6} -810000\,{y}^{4}{z}^{8} = 0.
\end{align*}%
}%
According to Maple, this curve is irreducible
over the algebraic closure of $\Q$.
\end{example}

\begin{prop}
The above curve of degree 12, together with the lines $x=y$ and $x=-y$,
is the branch curve of a double plane
whose smooth minimal model is a numerical Godeaux surface $S$
of Du Val type with $\Tors(S)=\Z/4\Z$.
\end{prop}

\begin{proof}
By Theorem \ref{settica},
it remains only to check that there exists a plane curve of degree 8
with the singularities described in that statement.
Again using Maple, we found that
such a curve actually exists and is the union of the sextic
\begin{align*}
& 900\,{y}^{6}-300\,{y}^{4}xz-1800\,{y}^{4}{z}^{2}-519\,{y}^{4}{x}^{2}+
900\,{y}^{2}{z}^{4}+106\,{y}^{2}{x}^{3}z+1597\,{y}^{2}{x}^{2}{z}^{2}+\\
& -365\,{y}^{2}{x}^{4}+300\,{y}^{2}x{z}^{3}+130\,{x}^{5}z-884\,{x}^{2}{z}
^{4}+299\,{x}^{4}{z}^{2}-364\,{x}^{3}{z}^{3} =0.
\end{align*}
and the lines $y=z$, $y=-z$.
\end{proof}

In \cite{We3}, Caryn Werner shows that all numerical Godeaux surfaces
of Campedelli type, whose branch curve is invariant under a linear involution
of the plane, have $\Tors(S)=\Z/4\Z$.
It would be interesting to see whether
a similar statement is still true for
numerical Godeaux surfaces of Du Val type.

The following example is computationally more complicated,
but it gives also instances of numerical Godeaux surfaces $S$
of Du Val type with $\Tors(S)=\Z/2\Z$.

\begin{example} \label{ex:deg11}
Let $[x,y,z]$ be homogeneous coordinates in $\pp^2$.
Consider the lines $r_1:x=0$, $r_2:x=y$, and the points
\[
q_0 = [0,0,1], q_1 = [0,1,0],  q_2 = [1,1,0],
q_3 = [1,0,0], q_4 = [-1,1,1], q_5 = [1,-2,1].
\]
We want to find a curve $B^\bullet$ which is reducible
in the line $r_3:y=0$, which passes through $q_0$ and $q_3$,
and in a curve $B_1^\bullet$ of degree 11
with the following singularities:
\begin{itemize}
\item two points of type $[4,4]$ in $q_1$, $q_2$,
with tangent lines $r_1$, $r_2$;
\item two quadruple points at $q_4$ and $q_5$;
\item two triple points at $q_0$, $q_3$;
\item a tacnode at a point $q_6=[t,0,1]$ on $r_3$,
with tacnodal tangent $r_3$.
\end{itemize}
Clearly, we want $q_6$ to be different from $q_0$, i.e.\ $t\ne0$.

Note that the virtual dimension of the linear system of curves
of degree 11 with these singularities is $-1$.
Thus it is reasonable to expect that for finitely many points on $r_3$
there is a curve with the required singularities.
Using Maple, we found a polynomial $p(t)$ of degree 15 in $t$
such that, if $t_0$ is a root of $p(t)$,
then there is a curve $B_1^\bullet$ with at least
the above singularities, with $q_6=[t_0,0,1]\ne q_0$.
Luckily, the polynomial $p(t)$ factors over $\Q$, and there is an irreducible
factor of degree 5 and another of degree 10.
Maple is able to work out the computations on the corresponding algebraic extensions
of $\Q$ and to verify that the corresponding curves $B_1^\bullet$
have the required, and not worse, singularities
and are irreducible over the algebraic closure of $\mathbb{Q}$.

Unfortunately, as we said, the computations are complicated and the results
are very cumbersome and it is not case to exhibit here the explicit equations,
which can be found at the following Web address:

\noindent\hfil {\tt http://www.mat.uniroma2.it/\~{ }calabri/duValeqs.pdf}

It is interesting to remark that the solutions relative to the degree 5 factor
have torsion $\Z/4\Z$, whereas the once of the factor of degree 10 have torsion $\Z/2\Z$.
This is verified by applying Corollary \ref{c:settica},
i.e.\ by checking the existence or not of the curve of degree 8
with the appropriate singularities, which in this case is the union of the line
$r_3$ and a curve of degree 7.

Finally, one can check, again by using Maple,
that the unique cubic plane curve $\bar\Delta$ as in Corollary \ref{c:notCamp}
is irreducible and does not pass through the irrelevant singularities of $B^\bullet$,
which is just a point in $B^\bullet_1\cap r_3$ different from $q_1, q_3, q_6$.
Therefore these surfaces cannot be of Campedelli type
by Corollary \ref{c:notCamp}.
\end{example}

Similar computations can be tried also in order to find irreducible
curves $B^\bullet$, giving $\Tors(S)=\Z/2\Z$.
Unfortunately, Maple is not able to carry out all the computations
in a reasonable time.
However we will see in a moment, in Corollary \ref{l:irr},
that there are equisingular deformations of
our curve which are irreducible.

\begin{rem} \label{r:moduli}
Using Maple, we verified that, in both the above examples,
imposing all the required singularities but $q_4$,
then the quadruple point $q_4$ can be chosen only in finitely many ways.
The verification is performed by checking that
there is no solution to the problem of finding the quadruple point $q_4$
on a given line, say $y-z=0$ in Example \ref{ex:Z4} the line,
and $x+y=0$ in Example \ref{ex:deg11}.

This agrees with the following fact.
The linear system of curves of degree 12 with the given singularities
at all the points but $q_4$ has virtual dimension 8,
which coincides with the actual dimension.
This can be verified using Maple.
Imposing a further fixed quadruple point $q_4$ is 10 more conditions.
However, by moving $q_4$ with two parameters in the plane,
we expect only eight more conditions.
The above computation shows that this naive expectation is actually right.
\end{rem}

We want to finish by giving information about the number
of parameters on which our constructions of numerical Godeaux
surfaces of Du Val type depend.

\begin{prop} \label{p:moduli}
Let $B^\bullet$ be the plane curve of degree 12
as in one of the two examples above.
Let $\V$ be a complete family of plane curves of degree 12
which is maximal under the property that
its general element is an equisingular deformation of $B^\bullet$.
Then $\dim \V=13$.
\end{prop}

\begin{proof}
It suffices to prove that $\dim(\V/PGL(3,\C))=5$.
Up to projective transformations,
we can fix distinct points $q_0$, $q_1$, $q_2$
and $q_3$, so that $r_1=\overline{q_0q_1}$ and $r_2=\overline{q_0q_2}$
are also fixed.
Let $\LL$ be the linear system of curves of degree 12
with quadruple points at $q_0$, $q_3$ and points of type $[4,4]$ at $q_i$,
where the tangent line is $r_i$, $i=1,2$.
Consider the following quadruplets $(B, p, q, \xi)\in \LL \times \pp^2 \times \pp^2 \times \Delta$,
where:
\begin{itemize}
\item  $\Delta\subset\pp^2\times{\pp^2}^\star$ is the incidence correspondence
so that $\xi \in \Delta$ is a pair $(x,r)$, with $x$ a point of the line $r$;
\item  $B \in \LL$ is irreducible and reduced;
\item  $q_0$, $\ldots$, $q_3$, $p$, $q$, $x$ are all distinct;
\item  $B$ has 4-uple points at $p$, $q$ and a $[3,3]$-point at $x$
with tangent direction $r$.
\end{itemize}
Let $\bar\V$ be an irreducible component of the closure of the set
of such quadruplets, containing the point $\beta=(B^\bullet, q_4,q_5,(q_6,r_3))$.
Note that $\dim\bar\V\ge5$.
In fact, $\bar V$ is defined in a neighbourhood of $\beta$
inside the variety $\LL \times \pp^2 \times \pp^2 \times \Delta$,
of dimension at least 37, by 32 equations.
Since the projection on the first factor from $\bar\V$ to $\LL$ is generically
finite to its image, it suffices to prove that $\dim\bar\V=5$.

Consider then the projection $\rho:\bar\V\to \pp^2\times\Delta$
to the last two factors.
Remark \ref{r:moduli} shows that $\rho$ is generically finite.
This concludes the proof.
\end{proof}

\begin{cor} \label{l:5-dim}
The components of the moduli space of numerical Godeaux surfaces
with an involution of Du Val type,
containing our two examples \ref{ex:Z4} and \ref{ex:deg11},
are 5-dimensional.
\qed
\end{cor}

Another interesting consequence is the following:

\begin{cor} \label{l:irr}
In case of Example \ref{ex:deg11},
the general element in $\bar\V$ corresponds to an irreducible curve.
\end{cor}

\begin{proof}
Since, as we saw in the proof of Proposition \ref{p:moduli},
the map $\rho$ is generically finite, then it is also surjective.
This means that the point $q_6$ can be moved out of the line
passing through $q_0$ and $q_3$,
hence, for the general member of $\bar\V$,
the line cannot be a component of $B^\bullet$.
\end{proof}

\bigskip

\small{\noindent Alberto Calabri,
Dipartimento di Matematica,
Universit\`a degli Studi di Bologna,
Piazza di Porta San Donato 5,  I-40126 Bologna, Italy,
%
{\it e-mail}: calabri@dm.unibo.it

\medskip
\noindent Ciro Ciliberto,
Dipartimento di Matematica,
Universit\`a degli Studi di Roma ``Tor Vergata'',
Via della Ricerca Scientifica, I-00133 Roma, Italy,
%
{\it e-mail}:  cilibert@mat.uniroma2.it

\medskip
\noindent Margarida Mendes Lopes,
Departamento de  Matem\'atica, Instituto Superior T\'ecnico,
Universidade T{\'e}cnica de Lisboa,
Av.~Rovisco Pais, 1049-001 Lisboa, Portugal, \\
%
{\it e-mail}: mmlopes@math.ist.utl.pt
}

\end{document}